\documentclass{siamonline0516}
\pdfoutput=1

\usepackage{amssymb,enumerate,amsfonts,mathrsfs}

\usepackage[retainorgcmds]{IEEEtrantools}
\usepackage{color}
\usepackage{lipsum}
\usepackage{soul}

\usepackage{enumitem}
\usepackage{enumerate}
\usepackage[utf8]{inputenc}

\theorembodyfont{\it}

\newtheorem{conjecture}{Conjecture}

\theorembodyfont{\rm}

\definecolor{cooperative}{rgb}{0.78,0.86,0.94}
\definecolor{competitive}{rgb}{0.68,0.72,0.78}
\usepackage{framed}
\usepackage{lipsum}
\usepackage{sidecap}

\def\dotquote{\begin{quote} }
\def\enddotquote{''\end{quote}}

\ifpdf
  \DeclareGraphicsExtensions{.eps,.pdf,.png}
\else
  \DeclareGraphicsExtensions{.eps}
\fi

\newcommand{\TheTitle}{A three-scale model of spatio-temporal bursting} 
\newcommand{\TheAuthors}{A. Franci and R. Sepulchre}

\headers{\TheTitle}{\TheAuthors}

\title{{\TheTitle}\thanks{Submitted to the editors 30 October 2015.
\funding{The research leading to these results has received funding from the European Research Council under the Advanced ERC Grant Agreement Switchlet n.670645 and from DGAPA-Universidad Nacional Autónoma de México under the PAPIIT Grant n.IA105816.}}}

\author{
  Alessio Franci\thanks{Department of Mathematics, Universidad Nacional Autónoma de México (UNAM), Ciudad de México, Mexico (\email{afranci@ciencias.unam.mx},\url{https://sites.google.com/site/francialessioac/}).}
  \and
  Rodolphe Sepulchre\thanks{Department of Engineering, University of Cambridge, Cambridge, UK (\email{r.sepulchre@eng.cam.ac.uk},\url{https://sites.google.com/site/rsepulchre/}).} 
}

\usepackage{amsopn}

\ifpdf
\hypersetup{
  pdftitle={\TheTitle},
  pdfauthor={\TheAuthors}
}
\fi

\begin{document}
\newcommand{\hy}{{\rm hy}}
\newcommand{\wc}{{\rm wcusp}}
\maketitle

% REQUIRED
\begin{abstract}
We study spatio-temporal bursting in a three-scale reaction diffusion equation organized by the winged cusp singularity. For large time-scale separation the model exhibits traveling bursts, whereas for large space-scale separation the model exhibits standing bursts. Both behaviors exhibit a common singular skeleton, whose geometry is fully determined by persistent bifurcation diagrams of the winged cusp. The modulation of spatio-temporal bursting in such a model naturally translates into paths in the universal unfolding of the winged-cusp.
\end{abstract}

% REQUIRED

% REQUIRED

\section{Introduction}

In the recent paper \cite{Franci2014}, we proposed that temporal bursting can be conveniently studied
in a three-time-scale model with a single scalar {\it winged-cusp} nonlinearity, a basic
organizing center in the language of  singularity theory. Because of its particular structure,
the model is readily interpreted as a three-scale generalization of the celebrated FitzHugh-Nagumo
model, a two-time-scale model similarly organized by a single scalar {\it hysteresis} nonlinearity.

With enough time-scale separation, the (temporal) patterns exhibited by the trajectories of those
models are robust and modulable in the parameter space, because they are in
one-to-one correspondence with the persistent bifurcation diagrams of the
universal unfolding of the hysteresis and winged-cusp, respectively.

We showed that this feature is of prime interest in neurophysiology, for instance, because
it provides insight in the physiological parameters that can regulate the continuous transition
between distinct temporal patterns as observed in experimental neuronal recordings.

The present paper pursues this investigation by adding diffusion in the three-scale temporal model
and exploring the resulting three-scale spatio-temporal patterns: traveling bursts and
standing bursts.

The close relationship of the model with its two-dimensional counterpart is of great help
in the analysis because we extensively rely on the existing theory for two-scale models with a hysteresis nonlinearity.
We study the existence of traveling bursts and standing bursts in our three-scale model
by mimicking the analysis of traveling pulses and standing pulses in earlier two-scale models. This
geometric analysis exploits the singular limit to construct the skeleton of the attractor
and relies on Fenichel theory to show the persistence of the pattern away from the singular
limit.

The results that we present are in close analogy both with the temporal results of the
three-dimensional model \cite{Franci2014} and with the spatiotemporal results in two-scale models \cite{Carpenter1977,Jones1984,Dockery1992,Jones1995a,Jones1998}.
We prove existence results when they follow from the existing Fenichel theory and discuss
sensible predictions in situations that require an extension of the existing
geometric theory. The conceptual message of the paper is that singularity theory provides
a principled methodology for a geometric study of spatiotemporal patterns that can be
robustly modulated in the presence of two or more well separated temporal and/or spatial scales.
The existence of specific patterns is proven for parameter sets sufficiently close
to a singular configuation (in the present situation, a pitchfork singularity in the universal unfolding of the winged cusp)
by enforcing the same geometric conditions as the classical two scale analysis.

Multi-scale patterns seem particularly relevant in biology. Traveling bursts play a major role in neuronal
dynamics, both during development \cite{Maeda1995} and in the adult \cite{Golomb1997}.
Morphogenetic fields with multiple characteristic spatial scales are hypothesized, for instance,
in limb development \cite{Newman1979}. More generally, we suggest that singularity theory provides a novel
perspective to model phenomena involving multi-scale communication via traveling waves and multi-scale pattern formation.

Mathematically, three components reaction diffusion equations were studied in \cite{Doelman2009}
but the organizing singularity was the hysteresis and the objects under analysis were single pulses, possibly interacting, rather than geometric bursts. Traveling bursts were also observed in a biophysical neuronal network model in \cite[Fig. 12.3]{Ermentrout2010}.
Averaging over the fast time-scale was used to reduce existence and stability analysis to the classical traveling pulse analysis.

The paper is organized as follows. We introduce the three-scale reaction-diffusion model under analysis in Section \ref{SEC: 3 scale model}. The system equations are prepared to be studied geometrically in Section \ref{SEC: methodology}, where we also recall some classical geometric tools for the analysis of traveling and standing pulses. Our presentation heavily relies on the ideas in \cite{Jones1995a,Jones1998}. We develop the fine (fast/short-range), medium (slow/long-range), and gross (ultraslow/ultra-long-range) scale analysis in Sections \ref{SEC: fine scale}, \ref{SEC: medium scale}, and \ref{SEC: gross scale}, respectively.  Section \ref{SEC: fine scale} revisits the existence of traveling and standing fronts around the winged-cusp nonlinearity. Sections \ref{SEC: medium scale} proves the co-existence of a homogenous resting state with a periodic wave-train or an infinite periodic pattern. Finally, the existence of traveling bursts is proved in Section \ref{SSECTravBurst} and the existence of standing bursts is conjectured (together with a detailed singular construction) in Section \ref{SSEC: standing bursts}. Based on this geometric analysis, and in full analogy with \cite{Franci2014},  Section \ref{SEC: modulation} discusses the relevance of the proposed model to study the modulation of spatio-temporal bursting. Open questions raised by the results of this paper are discussed theoretically and numerically in Section \ref{SEC: Discussion}.
 
\section{A three-scale reaction diffusion model}
\label{SEC: 3 scale model}

The paper studies spatio-temporal bursting in the reaction diffusion model
\begin{IEEEeqnarray}{rCl}\label{EQWcRd3D_temp}
\IEEEyesnumber
\tau_u u_t&=&D_u u_{xx}+g_{\wc}(u,\lambda+w,\alpha+z,\beta,\gamma),\IEEEyessubnumber\\
\tau_w w_t&=&D_w w_{xx} + u-w.\IEEEyessubnumber\\
\tau_z z_t&=&D_z z_{xx} + u-z\IEEEyessubnumber
\end{IEEEeqnarray}
where $\tau_u,\tau_w,\tau_z>0$, $D_u>0$, and $D_w,D_z\geq0$. The function $$g_{\wc}(u,\lambda,\alpha,\beta,\gamma)=-u^3-\lambda^2-\alpha+\beta u +\gamma u\lambda$$ is a universal unfolding of the winged cusp $-u^3-\lambda^2$. We regard (\ref{EQWcRd3D_temp}) as a one-component reaction-diffusion model organized by the winged-cusp singularity and with (spatio-temporal) adaptation variables $w$ and $z$.
%The specie $w$ models adaptation of the bifurcation parameter $\lambda$. The specie $z$ models adaptation of the unfolding parameter $\alpha$.
%The specie $u$ is an autocatalytic specie. Because $\frac{\partial g_{\wc}}{\partial w}$ is positive for small $w$ and negative for large $w$, the specie $w$ provides a mixed feedback on $u$. The specie provides a purely negative feedbac on $u$.\\

The fundamental hypothesis of model (\ref{EQWcRd3D_temp}) is a hierarchy of scales between the ``master" variable $u$ and the adaptation variables $w$ and $z$. This hierarchy of scales is consistent with the hierarchy between the state variable ($u$), the bifurcation parameter ($\lambda$), and unfolding parameters ($\alpha,\beta,\gamma$) in singularity theory applied to bifurcation problems \cite{Golubitsky1985}. The variable $w$ accounts for variations of the bifurcation parameter while the variable $z$ accounts for variations of the unfolding parameter $\alpha$. The fine-grain dynamics is organized by a single nonlinearity, with distinct behaviors depending on the bifurcation parameter and unfolding parameters. The medium-grain dynamics model (linear) adaptation of the bifurcation parameter to the quasi steady-state of the fine-grain behavior. The gross-grain dynamics models (linear) adaptation of unfolding parameters to the quasi steady-state behavior of the fine and medium grains.

We study separately the role of a time-scale separation
\begin{equation}\label{EQ:time-scale separation def}
0<\frac{\tau_u}{\tau_z}=:\varepsilon_{us}\ll\frac{\tau_u}{\tau_w}=:\varepsilon_s\ll 1,
\end{equation}
which leads to {\it traveling} burst waves, and the role of space-scale separation
\begin{equation}\label{EQ:spacescale separation def}
\sqrt{\frac{D_z}{D_u}}=:\delta_{ul}^{-1}\gg\sqrt{\frac{D_w}{D_u}}=:\delta_l^{-1}\gg 1,
\end{equation}
which leads to {\it standing} burst waves.

The proposed model extends earlier work in two directions:
\begin{itemize}
\item in the absence of diffusion, model (\ref{EQWcRd3D_temp}) reduces to the bursting model analyzed in \cite{Franci2014}. The results in \cite{Franci2014} underlie the importance of a persistent bifurcation diagram in the universal unfolding of the winged cusp, the {\it mirrored hysteresis} (class 2 of the persistent bifurcation diagram of the winged cusp listed in \cite[page 208]{Golubitsky1985})), as a key nonlinearity for obtaining robust three-time-scale bursting as ultra-slow adaptation around a slow-fast phase portrait with coexistence of a stable fixed point and a stable limit cycle \cite{Franci2012}.
\item when the winged cusp nonlinearity is replaced by an universal unfolding of the hysteresis singularity
$$g_{hy}(u,\lambda,\beta)=-x^3+\lambda+\beta x, $$
then (\ref{EQWcRd3D_temp}) reduces to a FitzHugh-Nagumo-type model, which is the fundamental two-scale model of excitability \cite{FitzHugh1961}, exhibiting traveling pulses under a time-scale separation \cite{Carpenter1977}\cite{Jones1984} and standing pulses under a space-scale separation \cite{Dockery1992}\cite{Jones1998}. 
\end{itemize}

\section{Methodology}
\label{SEC: methodology}

\subsection{Reduction to singularly perturbed ODEs}

Using the traveling wave anszatz $(u,w,z)\left(\frac{x}{\sqrt{D_u}}+c\frac{t}{\tau_u}\right)$, traveling and standing waves of (\ref{EQWcRd3D_temp}) satisfy the ordinary differential equation
\begin{IEEEeqnarray}{rCl}\label{EQWcRd3DTravODEtemp}
\IEEEyesnumber
u'&=&v_u \,,\IEEEyessubnumber\\
v_u'&=& cv_u-g_{\wc}(u,\lambda+w,\alpha+z,\beta,\gamma)\,,\IEEEyessubnumber\\
w'&=&v_w\,,\IEEEyessubnumber\\
v_w'&=& c \frac{\delta_l^2}{\varepsilon_s} v_w + \delta_l^2(w-u)\,,\IEEEyessubnumber\\
z'&=&v_z\,,\IEEEyessubnumber\\
v_z'&=& c \frac{\delta_{ul}^2}{\varepsilon_{us}} v_z + \delta_{ul}^2(z-u)\,,\IEEEyessubnumber
\end{IEEEeqnarray}

For $D_w,D_z\neq0$, we set $c=0$ and use the change of variables $v_{w}\mapsto\frac{\tilde v_{w}}{\delta_l},\ v_{z}\mapsto\frac{\tilde v_{z}}{\delta_{ul}}$ to study the existence of standing waves satisfying the ordinary differential equation
\begin{IEEEeqnarray}{rCl}\label{EQWcRd3DStandODE}
\IEEEyesnumber
u'&=&v_u\,,\IEEEyessubnumber\\
v_u'&=&-g_{\wc}(u,\lambda+w,\alpha+z,\beta,\gamma)\,,\IEEEyessubnumber\\
w'&=&\delta_l v_w \,,\IEEEyessubnumber\\
v_w'&=&\delta_l (  w - u)\,,\IEEEyessubnumber\\
z'&=&\delta_{ul} v_z \,,\IEEEyessubnumber\\
v_z'&=&\delta_{ul}( z -u)\,.\IEEEyessubnumber
\end{IEEEeqnarray}

To study traveling waves, we start by multiplying both sides of (\ref{EQWcRd3DTravODEtemp}d) and (\ref{EQWcRd3DTravODEtemp}f) by $\frac{\varepsilon_s}{c\delta_l^2}$ and $\frac{\varepsilon_{us}}{c\delta_{ul}^2}$, respectively,
\begin{IEEEeqnarray*}{rCl}
u'&=&v_u \,,\\
v_u'&=& cv_u-g_{\wc}(u,\lambda+w,\alpha+z,\beta,\gamma)\,,\\
w'&=&v_w\,,\\
\frac{\varepsilon_s}{c\delta_l^2} v_w'&=& v_w + \frac{\varepsilon_s}{c}(w-u)\,,\\
z'&=&v_z\,,\\
\frac{\varepsilon_{us}}{c\delta_{ul}^2} v_z'&=& v_z + \frac{\varepsilon_{us}}{c}(z-u)\,.
\end{IEEEeqnarray*}
In the limit $\frac{\varepsilon_s}{c\delta_l^2},\frac{\varepsilon_{us}}{c\delta_{ul}^2}\to0$, we obtain
\begin{IEEEeqnarray}{rCl}\label{EQWcRd3DTravODE}
\IEEEyesnumber
u'&=&v_u \,,\IEEEyessubnumber\\
v_u'&=& cv_u-g_{\wc}(u,\lambda+w,\alpha+z,\beta,\gamma)\,,\IEEEyessubnumber\\
w'&=&\frac{\varepsilon_s}{c} (u- w)\,,\IEEEyessubnumber\\
z'&=&\frac{\varepsilon_{us}}{c}( u-z)\,,\IEEEyessubnumber
%c'&=&0\,.\IEEEyessubnumber
\end{IEEEeqnarray}
Note that if $D_w=D_z=0$, that is, the two adaptation variable do not diffuse, then $\frac{\varepsilon_s}{c\delta_l^2},\frac{\varepsilon_{us}}{c\delta_{ul}^2}=0$,  because in this case $\delta_l=\delta_{ul}=+\infty$. The persistence of the limit  $\frac{\varepsilon_s}{c\delta_l^2},\frac{\varepsilon_{us}}{c\delta_{ul}^2}\to0$ for $D_w,D_z\neq 0$ will not be addressed in the present paper but we verify it numerically in Section~\ref{SEC: numerics}.

\subsection{The geometry of traveling and standing pulses}

For future reference, we briefly recall here the geometric construction of traveling and standing pulses in two-scale reaction diffusion equations. Our exposition is based on \cite[Sections 4.2, 4.5, and 5.3]{Jones1995a} for the traveling case and \cite[Sections 4 and 5]{Jones1998},\cite[Section 2]{Dockery1992} for the standing case.

Let $g_{hy}(u,\lambda,\beta)=-u^3+\beta u+\lambda$ be an universal unfolding of the hysteresis singularity. The two-scale reaction diffusion model
\begin{IEEEeqnarray}{rCl}\label{EQHyRd2D}
\IEEEyesnumber
\tau_u u_t&=&D_u u_{xx}+g_{\hy}(u,\lambda+w,\beta),\IEEEyessubnumber\\
\tau_w w_t&=&D_w w_{xx} + u-w,\IEEEyessubnumber
\end{IEEEeqnarray}
exhibits traveling pulses under time-scale separation (\ref{EQ:time-scale separation def}) and standing pulses under space-scale separation (\ref{EQ:spacescale separation def}). The traveling pulse is constructed geometrically for $D_w=0$ using the traveling wave ansatz.  The standing pulse is constructed geometrically in the traveling wave ansatz with $c=0$. For large time-scale separation and large space-scale separation the two reductions lead to the singularly perturbed ordinary differential equations

\begin{figure}[h!]
\centering
\includegraphics[width=0.8\textwidth]{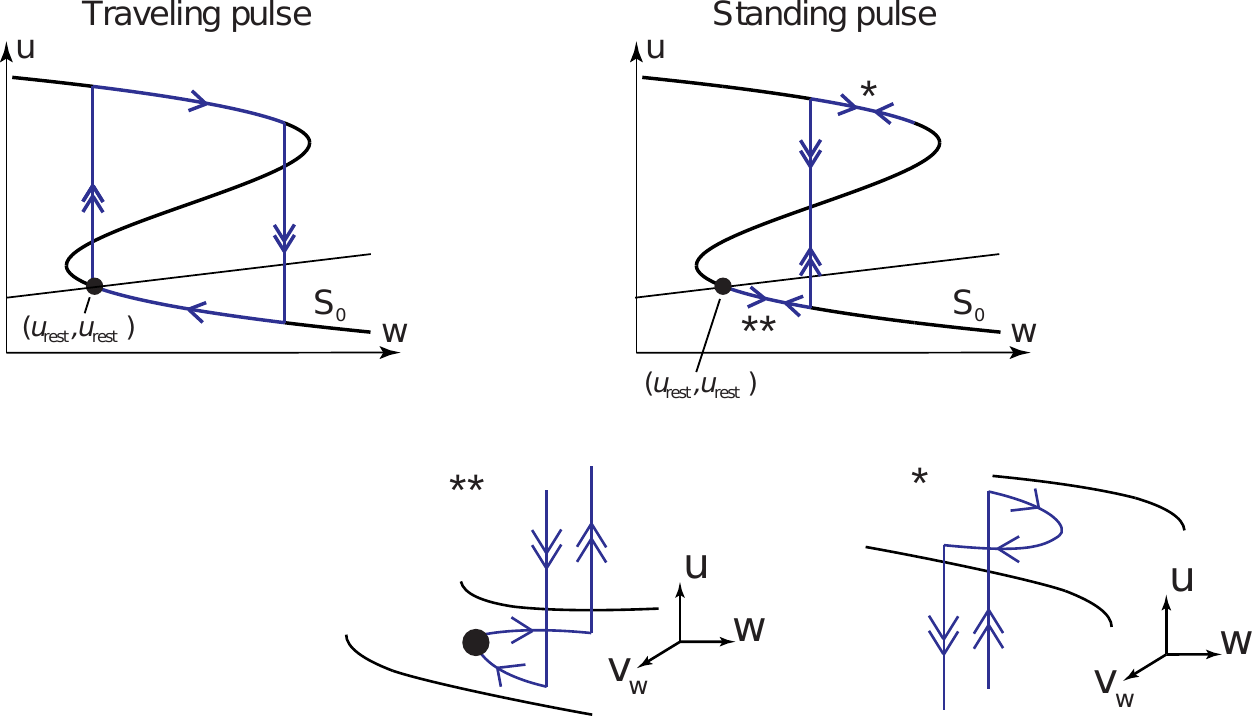}
\caption{Geometry of traveling and standing pulses in reaction-diffusion model with a cubic nonlinearity. Stable homogeneous resting states are depicted as black dots. The critical manifold $S_0$ is depicted as the thick black line. The thin black line is the nullcline of the slow variable $w$ (traveling case) or $v_w$ (standing case). Trajectories of the layer dynamics (\ref{EQWcRdFHNTravODElayer}) and (\ref{EQWcRdFHNStandODElayer}) are drawn with double arrows. Trajectories of the reduced dynamics (\ref{EQWcRdFHNTravODEreduced}) and (\ref{EQWcRdFHNStandODEreduced}) are drawn with single arrows.}\label{FIG:1FHN}
\end{figure}

\begin{center}
\begin{minipage}{0.4\linewidth}
\begin{IEEEeqnarray}{rCl}\label{EQWcRdFHNTravODE}
\IEEEyesnumber
u'&=&v_u \,,\IEEEyessubnumber\\
v_u'&=& cv_u-g_{\hy}(u,\lambda+w,\beta),\IEEEyessubnumber\\
w'&=&\frac{\varepsilon_s}{c} (u- w)\,,\IEEEyessubnumber\\
c'&=&0\,,\IEEEyessubnumber\\
&&\text{(traveling)}\nonumber\\
&&\nonumber
\end{IEEEeqnarray}
\end{minipage}
\begin{minipage}{0.5\linewidth}
\begin{IEEEeqnarray}{rCl}\label{EQWcRdFHNStandODE}
\IEEEyesnumber
u'&=&v_u\,,\IEEEyessubnumber\\
v_u'&=&-g_{\hy}(u,\lambda+w,\beta),\IEEEyessubnumber\\
w'&=&\delta_l v_w \,,\IEEEyessubnumber\\
v_w'&=&\delta_l (  w - u)\,.\IEEEyessubnumber\\
&&\text{(standing)}\nonumber\\
&&\nonumber
\end{IEEEeqnarray}
\end{minipage}
\end{center}
The dummy dynamics of the wave speed parameter $c$ in the traveling case is needed for the singular perturbation construction.
Under the same geometric conditions as \cite[Section 4.2]{Jones1995a} - traveling case - and \cite[Section 1]{Dockery1992} - standing case - the singular skeleton of (\ref{EQWcRdFHNTravODE}) and (\ref{EQWcRdFHNStandODE}) is sketched in Figure \ref{FIG:1FHN}. Geometrically, traveling and standing pulses share the same singular skeleton provided by their hysteretic {\it critical manifold}
%\begin{center}
%\begin{minipage}{0.49\linewidth}
%\begin{IEEEeqnarray}{l}\label{EQWcRdFHNTravODEcritmanifold}
%S_0=\nonumber\\
%\{(u,v_u,w,c):v_u=0,g_{\hy}(u,\lambda+w,\beta)=0\}\nonumber\\
%\quad\quad\text{(traveling)}\nonumber\\
%\nonumber
%\end{IEEEeqnarray}
%\end{minipage}
%\begin{minipage}{0.49\linewidth}
%\begin{IEEEeqnarray}{l}\label{EQWcRdFHNStandODEcritmanifold}
%S_0=\nonumber\\
%\{(u,v_u,w,v_w):v_u=v_w=0,g_{\hy}(u,\lambda+w,\beta)=0\}\nonumber\\
%\quad\quad\text{(standing)}\nonumber\\
%\nonumber
%\end{IEEEeqnarray}
%\end{minipage}
%\end{center}

\begin{IEEEeqnarray}{rCl}%\label{EQWcRdFHNTravODEcritmanifold}
S_0&=&\{(u,v_u,w,c):v_u=0,g_{\hy}(u,\lambda+w,\beta)=0\}\nonumber\\
&&\text{(traveling)}\nonumber\\
\nonumber\\
S_0&=&\{(u,v_u,w,v_w):v_u=v_w=0,g_{\hy}(u,\lambda+w,\beta)=0\}\nonumber\\
&&\text{(standing)}\nonumber
\end{IEEEeqnarray}

Both traveling and standing pulses are constructed as homoclinic orbits obtained as perturbations of the singular homoclinic orbits sketeched in Figure \ref{FIG:1FHN}.

The singular homoclinic orbit associated to the {\it traveling} pulse begins with a first jump from steady-state to the excited along the (fast) layer dynamics 
\begin{IEEEeqnarray}{rCl}\label{EQWcRdFHNTravODElayer}
\IEEEyesnumber
u'&=&v_u \,,\IEEEyessubnumber\\
v_u'&=& cv_u-g_{\hy}(u,\lambda+w,\beta)\,,\IEEEyessubnumber\\
w'&=&0\,,\IEEEyessubnumber\\
c'&=&0\,,\IEEEyessubnumber
\end{IEEEeqnarray}
This jump is a heteroclinic orbit of (\ref{EQWcRdFHNTravODElayer}) constructed as the transverse intersection in the 4-dimensional $(u,v_u,w,c)$-space of suitable 2-dimensional stable and 3-dimensional unstable manifolds (see \cite[Section 4.5]{Jones1995a} for details). Because $c$ is a parameter, this transverse intersection lies in the same $c$-slice, that is, for $c=c^*\neq0$. The phase portrait of the $(u,v_u)$ subsystem at the heteroclinic connection is the same as Figure~\ref{FIG:1}B$^{**}$. The slow motion inside the $1$-dimensional $c=c^*$-slice of the upper branch of the critical manifold is ruled by the reduced dynamics
\begin{IEEEeqnarray}{rCl}\label{EQWcRdFHNTravODEreduced}
\IEEEyesnumber
0&=&v_u \,,\IEEEyessubnumber\\
0&=& cv_u-g_{\hy}(u,\lambda+w,\beta)\,,\IEEEyessubnumber\\
w'&=&\frac{1}{c} (u- w)\,,\IEEEyessubnumber\\
c'&=&0,\quad c=c^*\,.\IEEEyessubnumber
\end{IEEEeqnarray}
The slow flow continues until the singular trajectory reaches the base point of a downward heteroclinic jump along the layer dynamics (\ref{EQWcRdFHNTravODElayer}). This heteroclinic orbit is again constructed as the transverse intersection of suitable stable and unstable manifolds (see \cite[Section 5.3]{Jones1995a} for details). The phase portrait of the $(u,v_u)$ subsystem at this heteroclinic connection is the same as Figure~\ref{FIG:1}B$^{*}$. The slow motion inside the $1$-dimensional $c=c^*$-slice of the lower branch of the critical manifold is ruled by the reduced dynamics (\ref{EQWcRdFHNTravODEreduced}) until the singular trajectory converges back to the resting point.

Due to the Hamiltonian nature of (\ref{EQWcRdFHNStandODE}), the {\it standing} pulse is a symmetric homoclinic orbit, that is, satisfying $(u(\xi),v_u(\xi),w(\xi),v_w(\xi))=(u(-\xi),-v_u(-\xi),w(-\xi),-v_w(-\xi))$. It therefore suffices to construct a singular trajectory originating at the steady state and intersecting transversely the subspace $\{v_u=v_w=0\}$. This singular homoclinic orbit begins with a slow portion along the reduced dynamics
\begin{IEEEeqnarray}{rCl}\label{EQWcRdFHNStandODEreduced}
\IEEEyesnumber
0&=&v_u\,,\IEEEyessubnumber\\
0&=&-g_{\hy}(u,\lambda+w,\beta)\,,\IEEEyessubnumber\\
w'&=&v_w \,,\IEEEyessubnumber\\
v_w'&=& w - u\,.\IEEEyessubnumber
\end{IEEEeqnarray}
Because $v_w'>0$ on the right of the steady state, $v_w$ and, therefore also $w$, increase until the trajectory reaches the base point of a heteroclinic trajectory along the layer dynamics
\begin{IEEEeqnarray}{rCl}\label{EQWcRdFHNStandODElayer}
\IEEEyesnumber
u'&=&v_u\,,\IEEEyessubnumber\\
v_u'&=&-g_{\hy}(u,\lambda+w,\beta)\,,\IEEEyessubnumber\\
w'&=&0,\IEEEyessubnumber\\
v_w'&=&0.\IEEEyessubnumber
\end{IEEEeqnarray}
Existence of this heteroclinic orbit easily follows from the Hamiltonian nature of the layer dynamics. The fact that it can be constructed as the transverse intersection in the four-dimensional $(u,v_u,w,v_w)$-space of suitable 2-dimensional stable and 3-dimensional unstable manifolds is proved by imposing the geometric conditions in \cite[Section 1]{Dockery1992}. At the heteroclinic jump the trajectory jumps to the upper branch of the critical manifold. There, the trajectory is carried transversely across $\{v_u=v_w=0\}$ by the reduced flow (\ref{EQWcRdFHNStandODEreduced}). Again, transversality, in particular the fact that the trajectory remains bounded away from the fold singularity of the critical manifold, is proved by imposing the geometric conditions in \cite[Section 1]{Dockery1992}. The other half of the trajectory is constructed by symmetry.  

Persistence of the singular homoclinic orbits associated to both the traveling and standing pulse follows by the Exchange Lemma \cite{Jones1994}. Roughly speaking, the Exchange Lemma allows to track, for $\varepsilon>0$, the invariant (stable and unstable) manifolds involved in the construction of the singular orbits and, in particular, to ensure that the same transverse intersections persist away from the singular limit and for $\varepsilon$ sufficiently small.

Application of the Exchange Lemma in the construction of the classical traveling and standing pulse requires that the homoclinic trajectory solely shadows normally hyperbolic part of the critical manifold. In the traveling case, this property is enforced by the fact that the homogeneous resting state is far from the fold singularity (note that the property may fail to hold when the resting state is close to the fold singularity or when $\varepsilon$ becomes too large \cite{Carter2015}). In the standing case it is enforced by the geometric assumptions in \cite[Section 1]{Dockery1992}. We will enforce the same property in the construction of our traveling and standing bursts by enforcing similar geometric conditions as the classical traveling and standing pulse.

\section{Fine-scale analysis: bistability and connecting fronts}
\label{SEC: fine scale}

For $\varepsilon_s=\varepsilon_{us}=0$ and $\delta_l=\delta_{ul}=0$, models (\ref{EQWcRd3DTravODE}) and (\ref{EQWcRd3DStandODE}) reduce to one-scale behaviors describing the fine-grain dynamics:
\begin{center}
\begin{minipage}{0.45\linewidth}
\begin{IEEEeqnarray}{rCl}\label{EQWcRd1DTravODE}
\IEEEyesnumber
u'&=&v_u \,,\IEEEyessubnumber\\
v_u'&=& cv_u-g_{\wc}(u,\tilde\lambda,\tilde\alpha,\beta,\gamma),\IEEEyessubnumber\\
&&\text{(traveling)}\nonumber\\
&&\nonumber
\end{IEEEeqnarray}
\end{minipage}
\begin{minipage}{0.45\linewidth}
\begin{IEEEeqnarray}{rCl}\label{EQWcRd1DStandODE}
\IEEEyesnumber
u'&=&v_u\,,\IEEEyessubnumber\\
v_u'&=&-g_{\wc}(u,\tilde\lambda,\tilde\alpha,\beta,\gamma),\IEEEyessubnumber\\
&&\text{(standing)}\nonumber\\
&&\nonumber
\end{IEEEeqnarray}
\end{minipage}
\end{center}
where $\tilde\lambda$, $\tilde\alpha$, $c$, $w$, and $z$ are now fixed parameters.

\begin{figure}[h!]
\centering
\includegraphics[width=0.9\textwidth]{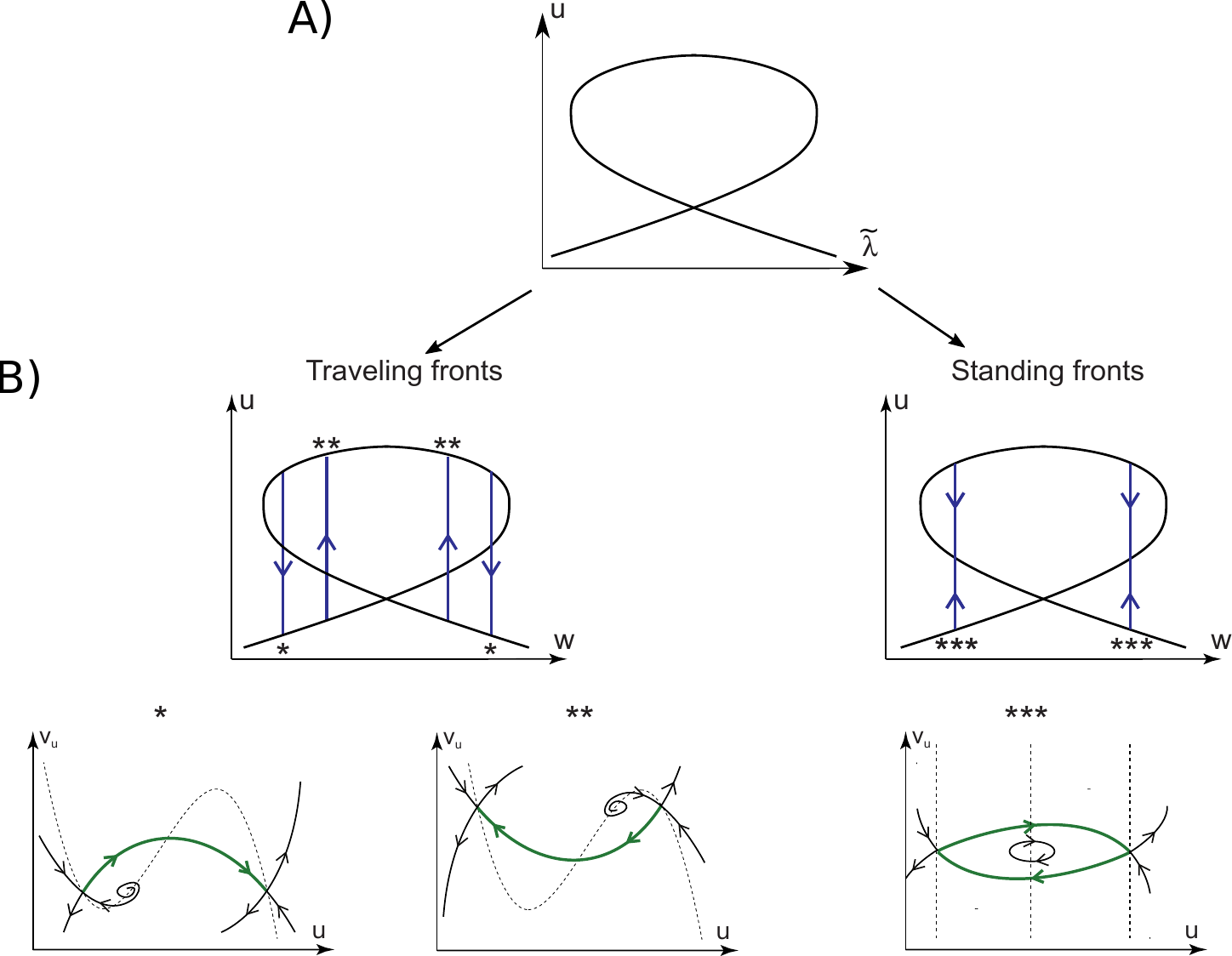}
\caption{The geometry of traveling and standing wave fronts around a transcritical singularity in the universal unfolding of the winged cusp. {\bf A.} Steady states of (\ref{EQWcRd1DTravODE}),(\ref{EQWcRd1DStandODE}) as a function of $\tilde\lambda$ for $\gamma=0$, $\beta>0$, and $\tilde\alpha=-2\left( \frac{\beta}{3}\right)^{2/3}$. {\bf B.} Top. Traveling and standing front of (\ref{EQWcRd1DTravODE}) and (\ref{EQWcRd1DStandODE}) respectively. Bottom. The $(u,v_u)$ phase portraits at the traveling and standing fronts.}\label{FIG:1}
\end{figure}

In both models, the steady-states
\begin{equation}\label{EQWcRd1D_ODESteady}
\{(u,v_u):\ v_u=0,\ g_{\wc}(u,\tilde\lambda,\tilde\alpha,\beta,\gamma)=0\}
\end{equation}
are determined by the universal unfolding of the winged cusp, as the bifurcation and unfolding parameters vary. For $\gamma=0$, $\beta>0$, and $\tilde\alpha=-2\left( \frac{\beta}{3}\right)^{2/3}$, the steady-state curve possesses a transcritical singularity for $\tilde\lambda=0$, where two mirror-symmetric hysteretic branches merge, as sketched in Figure~\ref{FIG:1}-A.

Away from the transcritical singularity, both hysteretic branches possess the same qualitative geometry as the classical cubic critical manifold in Figure \ref{FIG:1FHN}. In particular, they exhibit bistability between up and down homogeneous steady-states of the associated scalar reaction diffusion equation. Both in the traveling and standing case, there exist connecting heteroclinic orbits between the up and down steady-state branches (Figure~\ref{FIG:1}-B), which correspond to traveling or standing wave fronts of the associated reaction diffusion equation. In the traveling case, there are four $c$-dependent values of the bifurcation parameter $\tilde\lambda$ for which the model posses a heteroclinic orbit: two upward heteroclinic orbits and two downward heteroclinic orbits. In the standing case, there are two values of the bifurcation parameter $\tilde\lambda$ for which the model possesses both a downward and an upward heteroclinic orbit. The phase-portraits of the fast heteroclinic jumps are sketched in the insets.

Existence of these heteroclinic orbits follows exactly from the theory in \cite[Section 4.5]{Jones1995a} for the traveling case and in \cite[Section 4]{Jones1998} for the standing case. They all persist to parameter variations, in particular, to unfolding of the transcritical singularity. We omit here a detailed proof but the key derivations are recalled in Lemma \ref{LEM: sing hetero connection} for the traveling case and in Lemma \ref{LEM:stand front constr 1} for the standing case.

\section{Medium-scale analysis: bistability between homogeneous and periodic states}
\label{SEC: medium scale}

To model slow adaptation of the bifurcation parameter $\lambda$, we unfreeze the slow variable $w$ and study coexistence of a resting state (corresponding to a stable homogeneous resting state in the original PDE) and a limit cycle (corresponding to a periodic wave train or an infinite periodic patter in the original PDE) in the models
\begin{center}
\begin{minipage}{0.5\linewidth}
\begin{IEEEeqnarray}{rCl}\label{EQWcRd2DTravODE}
\IEEEyesnumber
u'&=&v_u\,,\IEEEyessubnumber\\
v_u'&=&cv_u-g_{\wc}(u,\lambda+w,\tilde\alpha,\beta,\gamma),\IEEEyessubnumber\\
w'&=&\frac{\varepsilon_s}{c}(u-w)\,,\IEEEyessubnumber\\
c'&=&0\,,\IEEEyessubnumber\\
&&\text{(traveling)}\nonumber\\
&&\nonumber
\end{IEEEeqnarray}
\end{minipage}
\begin{minipage}{0.45\linewidth}
\begin{IEEEeqnarray}{rCl}\label{EQWcRd2DStandODE}
\IEEEyesnumber
u'&=&v_u\,,\IEEEyessubnumber\\
v_u'&=&-g_{\wc}(u,\lambda+w,\tilde\alpha,\beta,\gamma),\IEEEyessubnumber\\
w'&=&\delta_l v_w\,,\IEEEyessubnumber\\
v_w'&=&\delta_{l} (w-u)\,,\IEEEyessubnumber\\
&&\text{(standing)}\nonumber\\
&&\nonumber
\end{IEEEeqnarray}
\end{minipage}
\end{center}
for $\varepsilon_s$ and $\delta_l$ sufficiently small, respectively. This bistability is the generalization of bistability between the homogeneous rest and excited states in the two-scale scenario.

For future reference, we make the following preliminary observation.  The homogeneous resting state equation of (\ref{EQWcRd2DTravODE}) and (\ref{EQWcRd2DStandODE}) is
\begin{equation}\label{EQHomSsEq}
F(u,\lambda,\tilde\alpha,\beta,\gamma):=-u^3-(\lambda+u)^2+\beta u-\gamma(\lambda+u)u-\tilde\alpha=0
\end{equation}
and is easily shown to be again a universal unfolding of the winged cusp around $u_{wcusp}:=\frac{1}{3}$, $\lambda_{wcusp}:=0$, $\alpha_{wcusp}:=-\frac{1}{27}$, $\beta_{wcusp}:=-\frac{1}{3}$, $\gamma_{wcusp}:=-2$.

\subsection{Traveling case}

\begin{figure}[h!]
\centering
\includegraphics[width=0.8\textwidth]{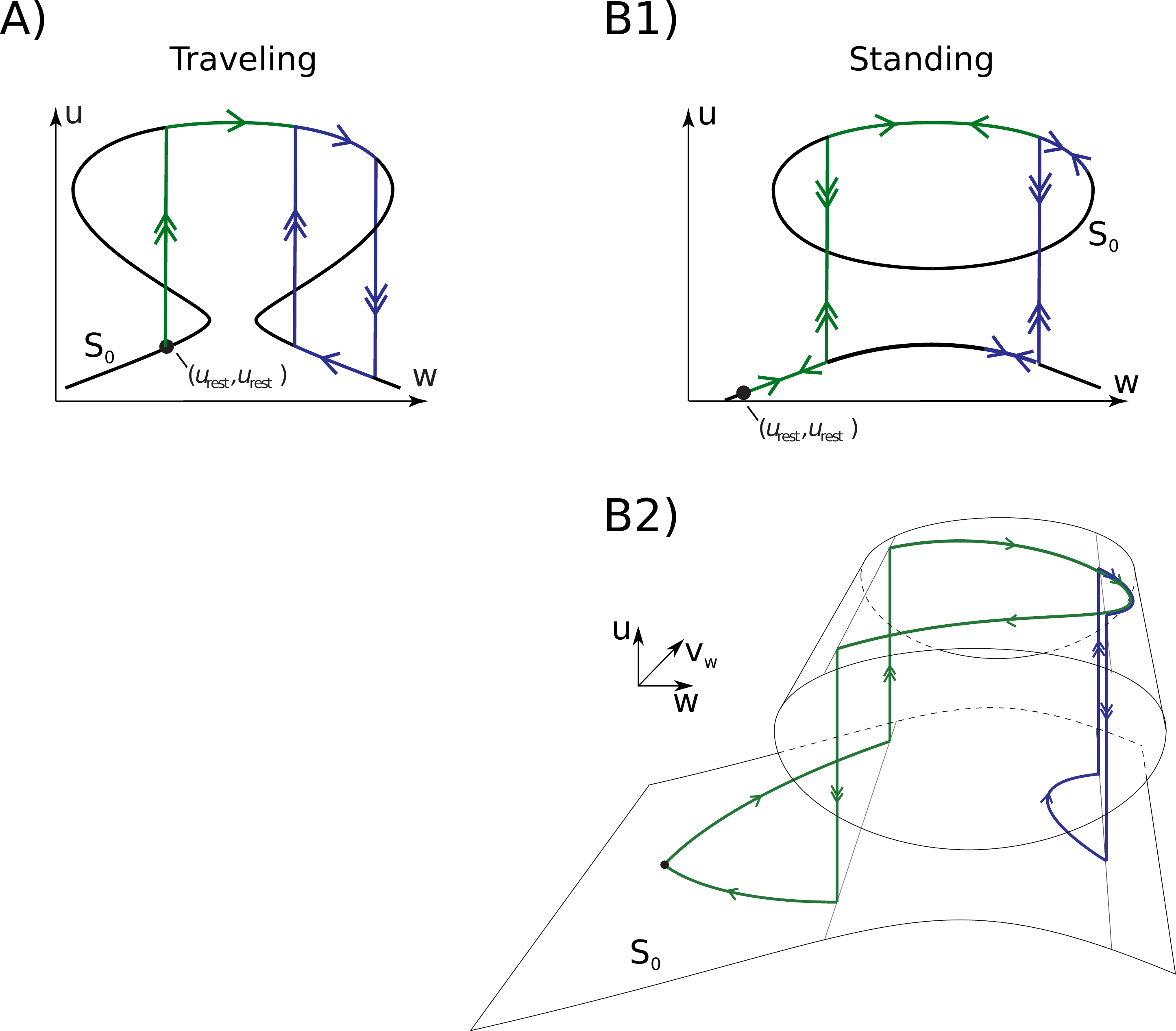}
\caption{{\bf A.} Bistability between a homogeneous steady state and a periodic wave-train in the singular limit and in the traveling wave ansatz. The two objects are connected by a traveling front. The homogeneous resting state is depicted as a black dot, the singular heteroclinic orbit corresponding to the traveling front as a green oriented line, and the singular periodic orbit corresponding to the periodic wavetrain as a blue oriented line. One arrow indicates slow portions along the reduced dynamics (\ref{EQWcRd2DTravODEslow}), two arrows fast portions along the layer dynamics (\ref{EQWcRd2DTravODEfast}).  {\bf B.} Bistability between a homogeneous steady state and an infinite periodic pattern in the singular limit and in the standing wave ansatz. There is also a (single) standing pulse solution that shadows the periodic pattern. The homogeneous resting state is depicted as a black dot, the singular periodic orbit corresponding to the infinite periodic pattern as a blue oriented line, and the singular homoclinic orbit corresponding to the standing pulse as a green oriented line. One arrow indicates slow portions along the reduced dynamics (\ref{EQWcRd2DStandODEslow}), two arrows fast portions along the layer dynamics (\ref{EQWcRd2DStandODEfast}).}\label{FIG3}
\end{figure}

The two-dimensional critical manifold of the singularly perturbed dynamics (\ref{EQWcRd2DTravODE}) is
\begin{equation}\label{EQWcRd2DTravODECritMan}
S_0=\{(u,v_u,w,c):\ v_u=0,\ g_{\wc}(u,\lambda+w,\tilde\alpha,\beta,\gamma)=0\}
\end{equation}
For $\gamma=0$, $\beta>0$, and $\tilde\alpha<-2\left( \frac{\beta}{3}\right)^{2/3}$ the transcritical singularity in Figure \ref{FIG:1}-A unfolds into the mirrored hysteresis persistent bifurcation diagram introduced in \cite[Section 2.3]{Franci2014} and organizing the singular phase portrait in Figure \ref{FIG3}-A. This algebraic curve is a template for rest-spike bistability in ordinary differential equations \cite[Theorem 2]{Franci2012}. Here we show that it is also a template to construct a wavefront between homogeneous resting and periodic wave-train.

The singular phase portrait of (\ref{EQWcRd2DTravODE}) in Fig \ref{FIG3}-A is constructed in  Lemma~\ref{LEM: sing hetero connection}. It exists for a precise value $c=c^*$ derived along the same line as \cite[Section 4.2]{Jones1995a}. It contains an equilibrium (black dot) corresponding to a stable homogeneous resting state of the associated reaction diffusion model, a singular periodic solution (blue line) corresponding to a traveling periodic wave train, and a heteroclinic\footnote{The term heteroclinic is used in the generalized sense of a trajectory connecting two $\omega$-limit sets, not necessarily two fixed points.} orbit (green line) corresponding to a traveling wave front connecting the stable homogeneous resting state to the periodic wave train. The singular periodic solution consists of two heteroclinic jumps of the layer dynamics (\ref{EQWcRd2DTravODEfast}) connected by trajectories of the reduced dynamics (\ref{EQWcRd2DTravODEslow}):
\begin{center}
\begin{minipage}{0.49\linewidth}
\begin{IEEEeqnarray}{rCl}\label{EQWcRd2DTravODEfast}
\IEEEyesnumber
u'&=&v_u\,,\IEEEyessubnumber\\
v_u'&=&cv_u-g_{\wc}(u,\lambda+w,\tilde\alpha,\beta,\gamma),\IEEEyessubnumber\\
w'&=&0\,,\IEEEyessubnumber\\
c'&=&0\IEEEyessubnumber\\
&&\nonumber
\end{IEEEeqnarray}
\end{minipage}
\begin{minipage}{0.49\linewidth}
\begin{IEEEeqnarray}{rCl}\label{EQWcRd2DTravODEslow}
\IEEEyesnumber
0&=&v_u\,,\IEEEyessubnumber\\
0&=&cv_u-g_{\wc}(u,\lambda+w,\tilde\alpha,\beta,\gamma),\IEEEyessubnumber\\
w'&=&\frac{1}{c}(u-w)\,,\IEEEyessubnumber\\
c'&=&0\,,\IEEEyessubnumber\\
&&\nonumber
\end{IEEEeqnarray}
\end{minipage}
\end{center}
The singular heteroclinic orbit consists of a heteroclinic jumps of the layer dynamics followed by a trajectory of the reduced dynamics.

The heteroclinic jumps are constructed in Lemma~\ref{LEM: sing hetero connection} as the transverse intersection in the 4-dimensional $(u,v_u,w,c)$-space of suitable 2-dimensional stable and 3-dimensional unstable manifolds of the layer dynamics (\ref{EQWcRd2DTravODEfast}), by enforcing the same geometric conditions as \cite[Section 4.2]{Jones1995a}. Like the standing pulse case, because $c$ is a parameter, this transverse intersection lies in the same $c$-slice, that is, for $c=c^*\neq0$. The fast phase portrait of the $(u,v_u)$ subsystem associated to upward (resp. downward) heteroclinic jumps is exactly the same as in Figure \ref{FIG:1}-B$^*$ (resp. Figure \ref{FIG:1}-B$^{**}$). The slow flow along the reduced dynamics (\ref{EQWcRd2DTravODEslow}) connects the landing point of successive heteroclinic orbits. This singular structure persists for $\varepsilon_s>0$, as proved in the following theorem. Its proof is provided in Section~\ref{SSEC: trav front proof}.

\vspace{2.5mm}
\begin{theorem}\label{THM:Wavefront}
There exist open sets of bifurcation ($\lambda$) and unfolding ($\alpha,\beta,\gamma$) parameters in a parametric neighborhood of the pitchfork singularity of (\ref{EQHomSsEq}) at $\beta=\frac{1}{3}$, $\lambda=\frac{1}{3}$, $\tilde\alpha=\frac{1}{27}-\frac{1}{9}$, $\gamma=0$ such that, for all parameters in those sets, there exists $c^*\neq0$ and $\bar\varepsilon_s>0$ such that, for each $0<\varepsilon_s<\bar\varepsilon_s$, model (\ref{EQWcRd2DTravODE}) possesses for a value $\bar c=c^*+\mathcal O(\varepsilon_s)$ a fixed point and a limit cycle both of saddle type and a heteroclinic orbit $h_{\varepsilon_s}$ converging in backward time to the fixed point and in forward time to the limit cycle.
%With an abuse of notation, we call $h_\varepsilon$ a heteroclinic orbit.
\end{theorem}
\vspace{2.5mm}

Theorem \ref{THM:Wavefront} generalizes bistability between a fixed point and a limit cycle in finite dimensional systems \cite[Theorem 2]{Franci2014} to bistability between a homogeneous resting state and a periodic wave-train in reaction-diffusion models organized by the winged cusp singularity. We note that the heteroclinic orbit of Theorem~\ref{THM:Wavefront} selects one out an infinity of periodic solutions for (\ref{EQWcRd2DTravODE}) by selecting a unique traveling speed.

 In the proof of Theorem~\ref{THM:Wavefront} we enforce the same standing conditions as for the construction of the classical FitzHugh-Nagumo pulse and apply the Exchange Lemma similarly. The singular skeleton is constructed along the same line as \cite[Sections 4.2, 4.5]{Jones1995a}. Existence of heteroclinic and periodic orbits is then proved by applying the Exchange Lemma along the same line as \cite[Sections 5.3]{Jones1995a}.

\subsection{Standing case}

The two-dimensional critical manifold of the singularly perturbed dynamics (\ref{EQWcRd2DStandODE}) is
\begin{equation}\label{EQWcRd2DStandODECritMan}
S_0=\{(u,v_u,w,v_w,c):\ v_u=v_w=0,\ g_{\wc}(u,\lambda+w,\tilde\alpha,\beta,\gamma)=0\}
\end{equation}
For $\gamma=0$, $\beta>0$, $\tilde\alpha>-2\left( \frac{\beta}{3}\right)^{2/3}$, and $\left|\tilde\alpha+2\left( \frac{\beta}{3}\right)^{2/3}\right|$ sufficiently small, the transcritical singularity in Figure \ref{FIG:1}-A unfolds into the persistent bifurcation diagram organizing the singular phase portrait in Figure~\ref{FIG3}-B1. This persistent bifurcation diagram is qualitatively different from the mirrored hysteresis in Figure~\ref{FIG3}-A but both diagrams belong to the unfolding of the same transcritical singularity in the universal unfolding of the winged cusp.

The singular phase portrait in Figure~\ref{FIG3}-B is constructed in Lemma~\ref{LEM:stand front constr 1}. It contains an equilibrium (black dot) corresponding to a stable homogeneous resting state of the associated reaction diffusion model, a singular periodic solution (blue line) corresponding to a periodic pattern, and a singular homoclinic orbit (green line) corresponding to standing pulse. The singular periodic and the homoclinic solution share a common arc: the standing pulse shadows the periodic pattern in its excited phase. The 3D projection of this singular phase-portrait onto the $(u,w,v_2)$ space is sketched in Figure \ref{FIG3}-B2. Both the periodic and homoclinic solutions consist of two heteroclinic jumps of the layer dynamics (\ref{EQWcRd2DStandODEfast}) connected by trajectories of the reduced dynamics (\ref{EQWcRd2DStandODEslow}):
\begin{center}
\begin{minipage}{0.49\linewidth}
\begin{IEEEeqnarray}{rCl}\label{EQWcRd2DStandODEfast}
\IEEEyesnumber
u'&=&v_u\,,\IEEEyessubnumber\\
v_u'&=&-g_{\wc}(u,\lambda+w,\tilde\alpha,\beta,\gamma)\,,\IEEEyessubnumber\\
w'&=&0\,,\IEEEyessubnumber\\
v_w'&=&0\,,\IEEEyessubnumber\\
&&\nonumber
\end{IEEEeqnarray}
\end{minipage}
\begin{minipage}{0.49\linewidth}
\begin{IEEEeqnarray}{rCl}\label{EQWcRd2DStandODEslow}
\IEEEyesnumber
0&=&v_u\,,\IEEEyessubnumber\\
0&=&-g_{\wc}(u,\lambda+w,\tilde\alpha,\beta,\gamma)\,,\IEEEyessubnumber\\
w'&=&v_w\,,\IEEEyessubnumber\\
v_w'&=&u-w\,.\IEEEyessubnumber\\
&&\nonumber
\end{IEEEeqnarray}
\end{minipage}
\end{center}
The fast phase portrait associated to the heteroclinic jumps is exactly the same as in Figure \ref{FIG:1} (case$^{***}$).

In analogy with \cite[Section 2]{Dockery1992} and \cite[Section 4]{Jones1998}, the first branch of the singular homoclinic orbit is provided by the slow 1-dimensional unstable manifold of the resting state inside the critical manifold, that is, according to the reduced dynamics (\ref{EQWcRd2DStandODEslow}). The slow flow continues until it reaches the fast heteroclinic trajectory of the layer dynamics (\ref{EQWcRd2DStandODEfast}). This trajectory, whose existence is easily proved by invoking the hamiltonian nature of (\ref{EQWcRd2DStandODEfast}), is constructed as the transverse intersection in the 4-dimensional $(u,v_u,w,v_w)$-space of suitable 2-dimensional unstable and 3-dimensional stable manifolds, much in the same way as \cite[Lemma 4.1]{Jones1998}. At the heteroclinic jump the trajectory jumps to the upper branch of the critical manifold. There, the trajectory is carried transversely across $\{v_u=v_w=0\}$ by the reduced flow. Avoidance of the critical manifold fold singularity is ensured by enforcing the same geometric conditions as  \cite[Section 1]{Dockery1992}. The second half of the singular homoclinic trajectory is constructed by symmetry of (\ref{EQWcRd2DStandODEfast}-\ref{EQWcRd2DStandODEslow}).

The singular periodic trajectory is constructed around the second (rightmost) fast heteroclinic connection, whose existence and transversality conditions are the same as the leftmost heteroclinic. The slow arc of the singular homoclinic trajectory lying in the upper branch of the critical manifold, on the right of the base points of the second fast heteroclinic connection provides the upper slow portion of the singular periodic trajectory. This is the shared arc. The fast heteroclinic connections connect this arc with a slow arc in the lower branch of the critical manifold.

This singular structure persists for $\delta_l>0$ as proved in the next Theorem. Its proof is provided in Section~\ref{SSEC: stand front proof}.

\vspace{2.5mm}
\begin{theorem}\label{THM:Standfront}
There exist open sets of bifurcation ($\lambda$) and unfolding ($\tilde\alpha,\beta,\gamma$) parameters
in a parametric neighborhood of the pitchfork singularity of (\ref{EQHomSsEq}) at $\beta=\frac{1}{3}$, $\lambda=\frac{1}{3}$, $\tilde\alpha=\frac{1}{27}-\frac{1}{9}$, $\gamma=0$
such that, for all parameters in those sets, there exists $\bar\delta_l>0$ such that, for all $0<\delta_l<\bar\delta_l$, model (\ref{EQWcRd2DStandODE}) possesses a homoclinic trajectory and a limit cycle that are $\mathcal O(\delta_l)$-close to each other together with their unstable manifolds in a neighborhood of the point $(u_{max},0,w_{max},0)$ where $w$ reaches its (unique) maximum along the singular homoclinic trajectory.
\end{theorem}
\vspace{2.5mm}

Theorem \ref{THM:Standfront} characterizes bistability between a homogeneous resting state and a periodic pattern in reaction diffusion models organized by the winged cusp singularity. We note that the homoclinic orbit singles out a privileged periodic pattern in model (\ref{EQWcRd2DStandODE}) out of an infinity of them.

In the proof of Theorem~\ref{THM:Standfront} we enforce the same conditions as the construction of the standing pulse and use the Exchange Lemma similarly. The singular skeleton is constructed by enforcing the same assumptions as \cite[Section 2]{Dockery1992}. The existence of the homoclinic and periodic orbits is proved along the same lines as \cite[Sections 4 and 5]{Jones1998}.

\section{Gross scale analysis: traveling and standing bursts}
\label{SEC: gross scale}

\subsection{Traveling bursts}
\label{SSECTravBurst}

\begin{figure}
\centering
\includegraphics[width=0.9\textwidth]{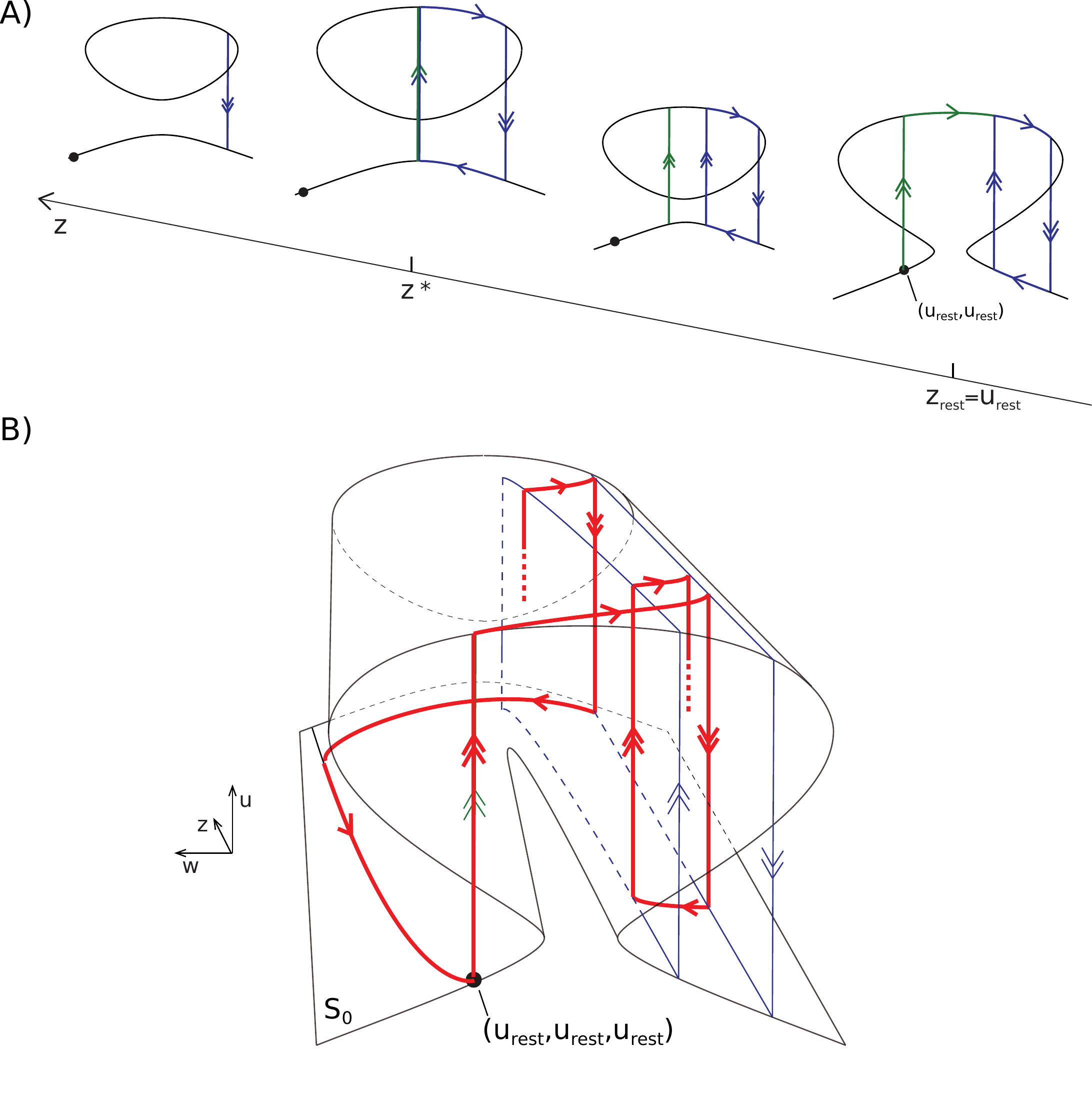}
\caption{Geometric construction of traveling bursts. {\bf A.} Deformation of the slow-fast phase portrait of (\ref{EQWcRd2DTravODE}) for different values of $\tilde\alpha=\alpha+z$, for fixed $\alpha$ and different values of $z$. Notation as Figure~\ref{FIG3}. {\bf B.} Projection of the singular homoclinic orbit corresponding to the traveling burst solution (in red) onto the $(u,w,z)$ space. One arrow indicates slow portions along the layer dynamics (\ref{EQWcRd3DTravODEreduced}), two arrows fast portions along the layer dynamics (\ref{EQWcRd3DTravODElayer}). Base and landing points of the fast heteroclinic connections lie along the blue lines.   } \label{FIG4}
\end{figure}

Very much like traveling pulses are homoclinic orbits of the singularly perturbed dynamics (\ref{EQWcRdFHNTravODE}), traveling bursts are constructed as homoclinic orbits of the singularly perturbed dynamics (\ref{EQWcRd3DTravODE}). We start by constructing the singular limit $\varepsilon_s\to0$ of this homoclinic orbit under the assumption that $\varepsilon_{us}=\tilde\varepsilon_{us}\varepsilon_{s}$, with $0<\tilde\varepsilon_{us}\ll 1$. The layer (\ref{EQWcRd3DTravODElayer}) and reduced (\ref{EQWcRd3DTravODEreduced}) dynamics then read
\begin{center}
\begin{minipage}{0.49\linewidth}
\begin{IEEEeqnarray}{rCl}\label{EQWcRd3DTravODElayer}
\IEEEyesnumber
u'&=&v\,,\IEEEyessubnumber\\
v'&=&cv-\nonumber\\
&&g_{\wc}(u,\lambda+w,\alpha+z,\beta,\gamma),\IEEEyessubnumber\\
w'&=&0\,,\IEEEyessubnumber\\
z'&=&0\,,\IEEEyessubnumber\\
c'&=&0\,,\IEEEyessubnumber\\
&&\nonumber
\end{IEEEeqnarray}
\end{minipage}
\begin{minipage}{0.49\linewidth}
\begin{IEEEeqnarray}{rCl}\label{EQWcRd3DTravODEreduced}
\IEEEyesnumber
0&=&v\,,\IEEEyessubnumber\\
0&=&cv-\nonumber\\
&&g_{\wc}(u,\lambda+w,\alpha+z,\beta,\gamma),\IEEEyessubnumber\\
w'&=&\frac{1}{c}(u-w)\,,\IEEEyessubnumber\\
z'&=&\frac{\tilde\varepsilon_{us}}{c}(u-z)\,,\IEEEyessubnumber\\
c'&=&0\,.\IEEEyessubnumber\\
&&\nonumber
\end{IEEEeqnarray}
\end{minipage}
\end{center}
The associated critical manifold is
\begin{equation}\label{EQ: trav burst crit man}
S_0:=\{(u,v,w,z,c):v=0,g_{\wc}(u,\lambda+w,\alpha+z,\beta,\gamma)=0\}.
\end{equation}
The overall geometry of the singular homoclinic orbit is illustrated in Figure \ref{FIG4}.

Figure \ref{FIG4}-A illustrates the deformation of the singular phase portrait of the medium scale dynamics (\ref{EQWcRd2DTravODE}) for different values of the unfolding parameter $\tilde\alpha=\alpha+z$ for fixed $\alpha$ and changing $z$. This deformation is proved in  Lemma~\ref{LEM: sing hetero connection 2}. The starting point (at $z=z_{rest}$) is determined by Theorem~\ref{THM:Wavefront}, which corresponds to the phase portrait in Figure \ref{FIG3}A. This situation is created by picking $(\lambda,\tilde\alpha,\beta,\gamma)$ and $c=c^*$ as in Theorem~\ref{THM:Wavefront}, and $\alpha=\tilde\alpha-u_{rest}$ so that $(u_{rest},u_{rest})$ is an equilibrium for $z=z_{rest}=u_{rest}$. Lemma \ref{LEM: sing hetero connection 2} also shows that, for $z\geq z_{rest}$, the fixed point and the limit cycle persist over a range $[z_{rest},z^*)$, but that the fast heteroclinic orbit from the lower to the upper branch of the critical manifold progressively slides toward the second upward heteroclinic orbit as $z\to z^*$. For $z>z^*$ only the fixed point and the downward heteroclinic orbit persist.

The singular homoclinic orbit is constructed by gluing together the $z$-slices of Figure~\ref{FIG4}A in the range $[z_{rest},z^*+\theta]$, for small $\theta>0$, as illustrated in Figure~\ref{FIG4}B. As for the traveling pulse, a first heteroclinic jump brings the orbit from rest to the excited state. As $z$ increases in the ultraslow scale, the slow flow governed by (\ref{EQWcRd3DTravODEreduced}) brings the trajectory to the base point of the next fast heteroclinic jump. The alternation of slow flows and fast heteroclinic jumps, along the family of singular limit cycles,  continues until $z>z^*$. The orbit then relaxes back to rest along the lower branch of the critical manifold.

The following theorem proves the existence of a homoclinic orbit that tracks the singular homoclinic orbit away from the singular limit. Its proof is provided in Section \ref{SSEC: trav burst proof}.

\vspace{2.5mm}
\begin{theorem}\label{THM: traveling burst}
There exist open sets of bifurcation ($\lambda$) and unfolding ($\alpha,\beta,\gamma$) parameters in a parametric neighborhood of the pitchfork singularity of (\ref{EQHomSsEq}) at $\beta=\frac{1}{3}$, $\lambda=\frac{1}{3}$, $\tilde\alpha=\frac{1}{27}-\frac{1}{9}$, $\gamma=0$ such that, for all parameters in those sets, there exists $c^*\neq0$ { such that, for almost all $\tilde\varepsilon_{us}>0$ sufficiently small and $c=c^*$, the singular limit (\ref{EQWcRd3DTravODElayer})-(\ref{EQWcRd3DTravODEreduced}) possesses a transverse singular homoclinic orbit as sketched in Figure~\ref{FIG4}. Furthermore, model (\ref{EQWcRd3DTravODE}) possesses for $\varepsilon_s>0$ sufficiently small a homoclinic orbit near the transverse singular homoclinic orbit.
%There also exists $\bar\varepsilon_s>0$ such that, for all $0<\varepsilon_s<\bar\varepsilon_s$, there exists $\bar c=c^*+\mathcal O(\varepsilon_s)$, such that model (\ref{EQWcRd3DTravODE}) possesses for $c=\bar c$ a homoclinic orbit obtained as a $\mathcal O(\varepsilon_s)$-perturbation of the transverse singular homoclinic solution.
}
\end{theorem}
\vspace{2.5mm}

The proof of Theorem~\ref{THM: traveling burst} is a direct application of the Exchange Lemma, more precisely, of the theorem in \cite[Section 4]{Jones1994}, which applies the Exchange Lemma to persistence of singular homoclinic orbits like the one in Figure~\ref{FIG4}. The presence of a third scale governed by $\varepsilon_{us}$ or, equivalently, by $\tilde\varepsilon_{us}$ is solely exploited to enforce the transversality conditions required by the application of the exchange lemma to (\ref{EQWcRdFHNTravODE}) via the singular limit (\ref{EQWcRd3DTravODElayer}-\ref{EQWcRd3DTravODEreduced}). {Genericity in $\tilde\varepsilon_{us}$ arises by imposing some of these transversality conditions.

Non-generic values of $\tilde\varepsilon_{us}$ include spike-adding bifurcations in which the number of spike per burst in the singular homoclinic solution changes: for $\tilde\varepsilon_{us}$ sufficiently large there is only one spike per burst because the slow flow brings $z$ above $z^*$ after just one jump; for decreasing $\tilde\varepsilon_{us}$ the number of jumps necessary to bring $z$ above $z^*$ increases monotonically and so does the number of spikes per burst. Spike-adding bifurcations and the rich dynamics they bring are well known in the purely temporal saddle-homoclinic bursting setting (see, e.g., \cite{Lee1999,Linaro2012}).}

\subsection{Standing bursts}
\label{SSEC: standing bursts}

\begin{figure}
\centering
\includegraphics[width=0.9\textwidth]{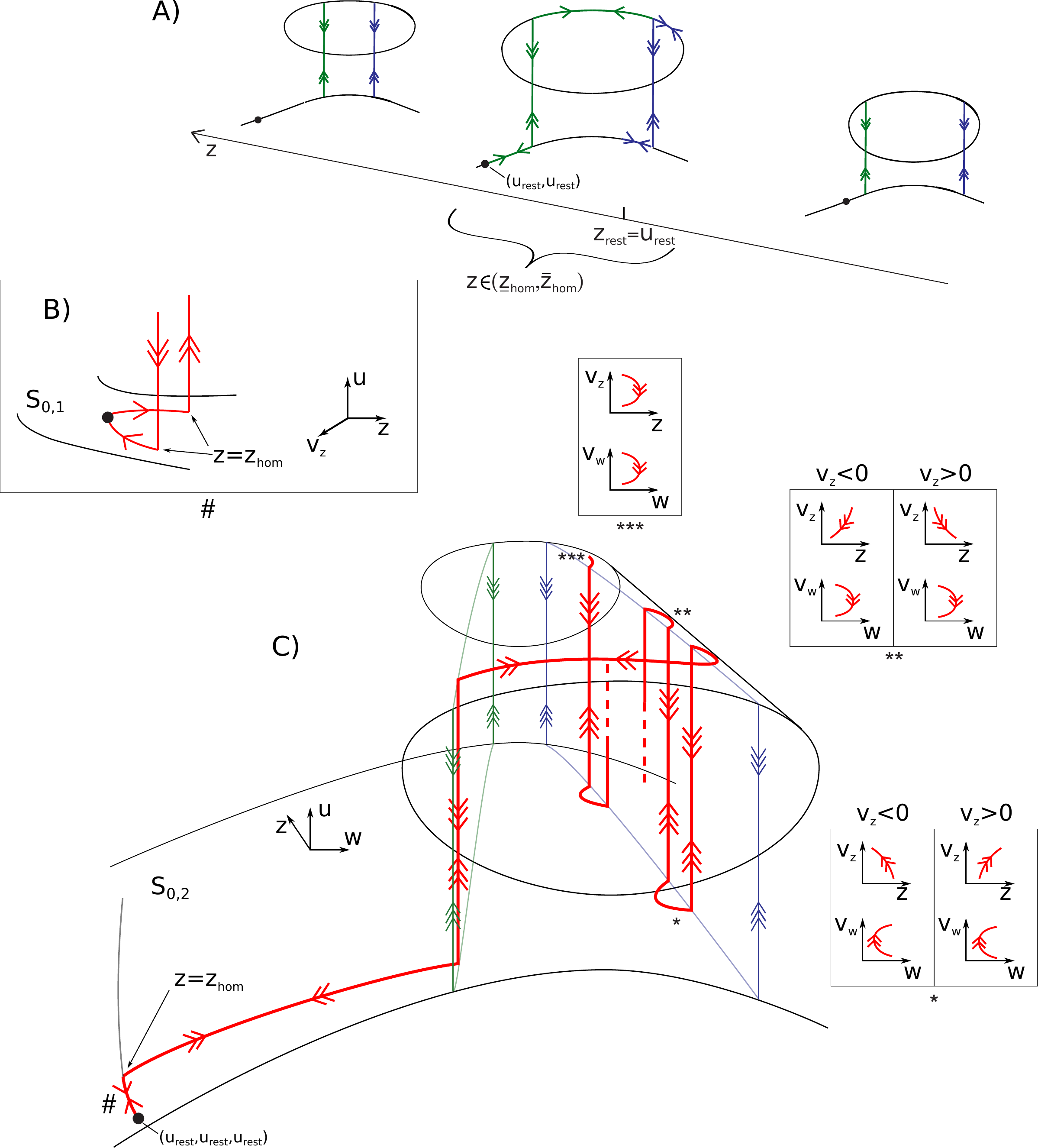}
\caption{Geometric construction of standing bursts. {\bf A.} Deformation of the slow-fast phase portrait of (\ref{EQWcRd2DStandODE}) for different values of $\tilde\alpha=\alpha+z$, for fixed $\alpha$ and different values of $z$. Notation as Figure~\ref{FIG3}. {\bf B.} and {\bf C.} Singular homoclinic orbit corresponding to the standing burst solution. One arrow distinguishes trajectories of the reduced dynamics (\ref{EQWcRd3DStandODEreducedA}), two arrows trajectories of the reduced dynamics (\ref{EQWcRd3DStandODEreducedB}), and three arrows trajectories of the layer dynamics (\ref{EQWcRd3DStandODElayerB}). Base and landing points of the fast heteroclinic connections lie along the blue lines. }\label{FIG5}
\end{figure}

In analogy with the standing pulses of (\ref{EQWcRdFHNStandODE}), the standing bursts of the singularly perturbed dynamics (\ref{EQWcRd3DStandODE})  are symmetric homoclinic orbits, that is, they satisfy
\begin{align*}
(u(\xi),v_u(\xi),w(\xi),v_w(\xi),z(\xi),v_z(\xi))=(u(-\xi),-v_u(-\xi),w(-\xi),-v_w(-\xi),z(-\xi),-v_z(-\xi)).
\end{align*}
The overall geometry of the associated symmetric {\it singular} homoclinic orbit is illustrated in Figure~\ref{FIG5}.

Figure~\ref{FIG5}A illustrates the deformation of the singular phase portrait of the medium scale dynamics (\ref{EQWcRd2DStandODE}) for different values of the unfolding parameter $\tilde\alpha=\alpha+z$, for fixed $\alpha$ and changing $z$.
%The existence of the fast heteroclinic connections is proved in Lemma~\ref{LEM:stand front constr 1}.
The existence of the range $(\underline z_{hom},\overline z_{hom})$, where a singular homoclinic orbit and a singular periodic orbit coexist, follows from Theorem~\ref{THM:Standfront}. In this range, the phase portrait is as in Figure \ref{FIG3}B. The condition $z_{rest}\in(\underline z_{hom},\overline z_{hom})$ is again imposed by picking $(\lambda,\tilde\alpha,\beta,\gamma)$ as in Theorem~\ref{THM:Standfront}, and $\alpha=\tilde\alpha-u_{rest}$ in such a way that $(u_{rest},u_{rest})$ is an equilibrium for $z=z_{rest}=u_{rest}$.
%Lemma~\ref{LEM:stand front constr 1} also shows that for $z^*=-\alpha>\overline z_{hom}$ no fast heteroclinic orbits persist.

The first portion of the singular symmetric homoclinic orbit is constructed in the singular limit $\delta_{ul}\to0$ and $0<\delta_l\ll1$, that is  $(u,v_u,w,v_w)$ are fast and $(z,v_z)$ are slow. The associated layer and reduced dynamics are
\begin{center}
\begin{minipage}{0.50\linewidth}
\begin{IEEEeqnarray}{rCl}\label{EQWcRd3DStandODElayerA}
\IEEEyesnumber
u'&=&v_u\,,\IEEEyessubnumber\\
v_u'&=&-g_{\wc}(u,\lambda+w,\alpha+z,\beta,\gamma),\IEEEyessubnumber\\
w'&=&\delta_l v_w\,,\IEEEyessubnumber\\
v_w'&=&\delta_l (w-u)\,,\IEEEyessubnumber\\
z'&=&0\,,\IEEEyessubnumber\\
v_z'&=&0\,,\IEEEyessubnumber\\
&&\nonumber
\end{IEEEeqnarray}
\end{minipage}
\begin{minipage}{0.49\linewidth}
\begin{IEEEeqnarray}{rCl}\label{EQWcRd3DStandODEreducedA}
\IEEEyesnumber
0&=&v_u\,,\IEEEyessubnumber\\
0&=&-g_{\wc}(u,\lambda+w,\alpha+z,\beta,\gamma),\IEEEyessubnumber\\
0&=&v_w\,,\IEEEyessubnumber\\
0&=&(w-u)\,,\IEEEyessubnumber\\
z'&=&v_z\,,\IEEEyessubnumber\\
v_z'&=&z-u-\bar z\,.\IEEEyessubnumber\\
&&\nonumber
\end{IEEEeqnarray}
\end{minipage}
\end{center}

In this singular limit, the homogeneous resting state belongs to the 2-dimensional critical manifold
\begin{equation*}
S_{0,1}:=\{(u,v_u,w,v_w,z,v_z):v_u=v_w=0,w=u,g_{\wc}(u,\lambda+w,\alpha+z,\beta,\gamma)=0\}
\end{equation*}
and, inside this manifold ({\it i.e.}, with respect to the reduced dynamics (\ref{EQWcRd3DStandODEreducedA})), it has one dimensional stable and unstable manifolds (Figure \ref{FIG5}B). The unstable manifold provides the first portion of the singular homoclinic orbit. Note that, in this singular limit, both $u$ and $w$ are at quasi steady state in this regime.

The reduced and layer dynamics (\ref{EQWcRd3DStandODElayerA}) and (\ref{EQWcRd3DStandODEreducedA}) are not appropriate to describe the oscillatory part of the pattern because there only $(u,v_u)$ are at their quasi steady-state (except during fast heteroclinic jumps). To construct the oscillatory part of the pattern we consider the singular limit $\delta_l\to 0$, $\delta_{ul}=\delta_l\tilde\delta_{ul}$ with $0<\tilde\delta_{ul}\ll1$, that is, $(u,v_u)$ are fast and $(w,v_w,z,v_z)$ are slow.
In this new singular limit, the layer and reduced dynamics are
\begin{center}
\begin{minipage}{0.49\linewidth}
\begin{IEEEeqnarray}{rCl}\label{EQWcRd3DStandODElayerB}
\IEEEyesnumber
u'&=&v_u\,,\IEEEyessubnumber\\
v_u'&=&-g_{\wc}(u,\lambda+w,\alpha+z,\beta,\gamma),\IEEEyessubnumber\\
w'&=&0\,,\IEEEyessubnumber\\
v_w'&=&0\,,\IEEEyessubnumber\\
z'&=&0\,,\IEEEyessubnumber\\
v_z'&=&0\,,\IEEEyessubnumber\\
&&\nonumber
\end{IEEEeqnarray}
\end{minipage}
\begin{minipage}{0.49\linewidth}
\begin{IEEEeqnarray}{rCl}\label{EQWcRd3DStandODEreducedB}
\IEEEyesnumber
0&=&v_u\,,\IEEEyessubnumber\\
0&=&-g_{\wc}(u,\lambda+w,\alpha+z,\beta,\gamma),\IEEEyessubnumber\\
w'&=&v_w\,,\IEEEyessubnumber\\
v_w'&=&w-u\,,\IEEEyessubnumber\\
z'&=&\tilde\delta_{ul} v_z\,,\IEEEyessubnumber\\
v_z'&=&\tilde\delta_{ul}(z-u-\bar z)\,.\IEEEyessubnumber\\
&&\nonumber
\end{IEEEeqnarray}
\end{minipage}
\end{center}
The four-dimensional critical manifold is the set
\begin{equation*}
S_{0,2}:=\{(u,v_u,w,v_w,z,v_z):v_u=0,g_{\wc}(u,\lambda+w,\alpha+z,\beta,\gamma)=0\}
\end{equation*}

The oscillatory part of the standing burst solution (Figure \ref{FIG5} C) is composed by trajectories of the reduced dynamics (\ref{EQWcRd3DStandODEreducedB}) connected by heteroclinic jumps of the layer dynamics  (\ref{EQWcRd3DStandODElayerB}). Projections of the slow orbits onto the $(w,v_w,z,v_z)$-space are drawn in the insets marked with stars.
As $z$ ultra-slowly increases, the trajectory shadows the family of singular periodic orbits. For $\tilde\delta_{ul}$ sufficiently small, this ensures that the slow portions also cross the subspace $\{v_w=0\}$ transversely.
The last slow orbit crosses transversely the subspace $\{v_w=v_z=0\}$, which makes the singular trajectory symmetric, in view of the fact that $v_u$ is identically zero on the critical manifold.

The quasi-steady state part of the standing burst, ruled by (\ref{EQWcRd3DStandODEreducedA}), connects to the oscillatory part of the standing burst, ruled by (\ref{EQWcRd3DStandODElayerB}), (\ref{EQWcRd3DStandODEreducedB}), at a value $z_{hom}$, where the singular trajectory leaves the quasi-steady branch. Let $W^u_{q.s.s}$ be the (three-dimensional) manifold obtained by the union of the (two-dimensional)  unstable manifolds of the fixed points of the layer dynamics (\ref{EQWcRd3DStandODElayerA}) lying on the quasi-steady state branch. Then $z_{hom}$ can be found by tracking $W^u_{q.s.s}$ along (\ref{EQWcRd3DStandODElayerB}), (\ref{EQWcRd3DStandODEreducedB}) and imposing its (zero-dimensional) transverse intersection with the (three-dimensional) $\{v_u=v_w=v_z=0\}$ subspace, much in the same way as \cite[Lemma 4.2]{Jones1998}.

The overall singular homoclinic orbit is constructed by patching two pieces corresponding to the two different singular limits. In contrast to traveling bursts, a proof of its persistence away from the singular limit(s) does not follow as an immediate application of the Exchange Lemma. However, due to transversality, its persistence is highly plausible, leading to the following conjecture.

\vspace{2.5mm}
\begin{conjecture}\label{THM: standing burst}
There exist open sets of bifurcation ($\lambda$) and unfolding ($\alpha,\beta,\gamma$) parameters in a parametric neighborhood of the pitchfork singularity of (\ref{EQHomSsEq}) at $\beta=\frac{1}{3}$, $\lambda=\frac{1}{3}$, $\tilde\alpha=\frac{1}{27}-\frac{1}{9}$, $\gamma=0$ such that, for all parameters in those sets {and almost all $\tilde\delta_{ul}>0$ sufficiently small, the double singular limit (\ref{EQWcRd3DStandODElayerA})-(\ref{EQWcRd3DStandODEreducedA}), (\ref{EQWcRd3DStandODElayerB})-(\ref{EQWcRd3DStandODEreducedB}) possesses a transverse singular homoclinic orbit as sketched in Figure \ref{FIG5}. Furthermore, model (\ref{EQWcRd3DStandODE}) possesses for $\delta_l>0$ sufficiently small a homoclinic orbit near the transverse singular homoclinic orbit. }
%there exists $\bar\delta_l>0$ such that, for all $0<\delta_l<\bar\delta_l$, there exists $\bar\delta_{ul}>0$ such that, for almost all $0<\delta_{ul}<\bar\delta_{ul}$, model (\ref{EQWcRd3DStandODE}) possesses a homoclinic orbit obtained as a $\mathcal O(\delta_{l})$-perturbation of the singular standing burst solution in Figure \ref{FIG5}.
\end{conjecture}

\section{Modulation of spatio-temporal bursting}
\label{SEC: modulation}
A main motivation in \cite{Franci2014} to study a model organized by the winged cusp singularity is that it provides a principled way to analyze the deformation of bursting patterns as parameter paths in the universal unfolding of the singularity.
The same geometric picture generalizes to spatio-temporal behaviors.

\begin{figure}
\centering
\includegraphics[width=0.6\textwidth]{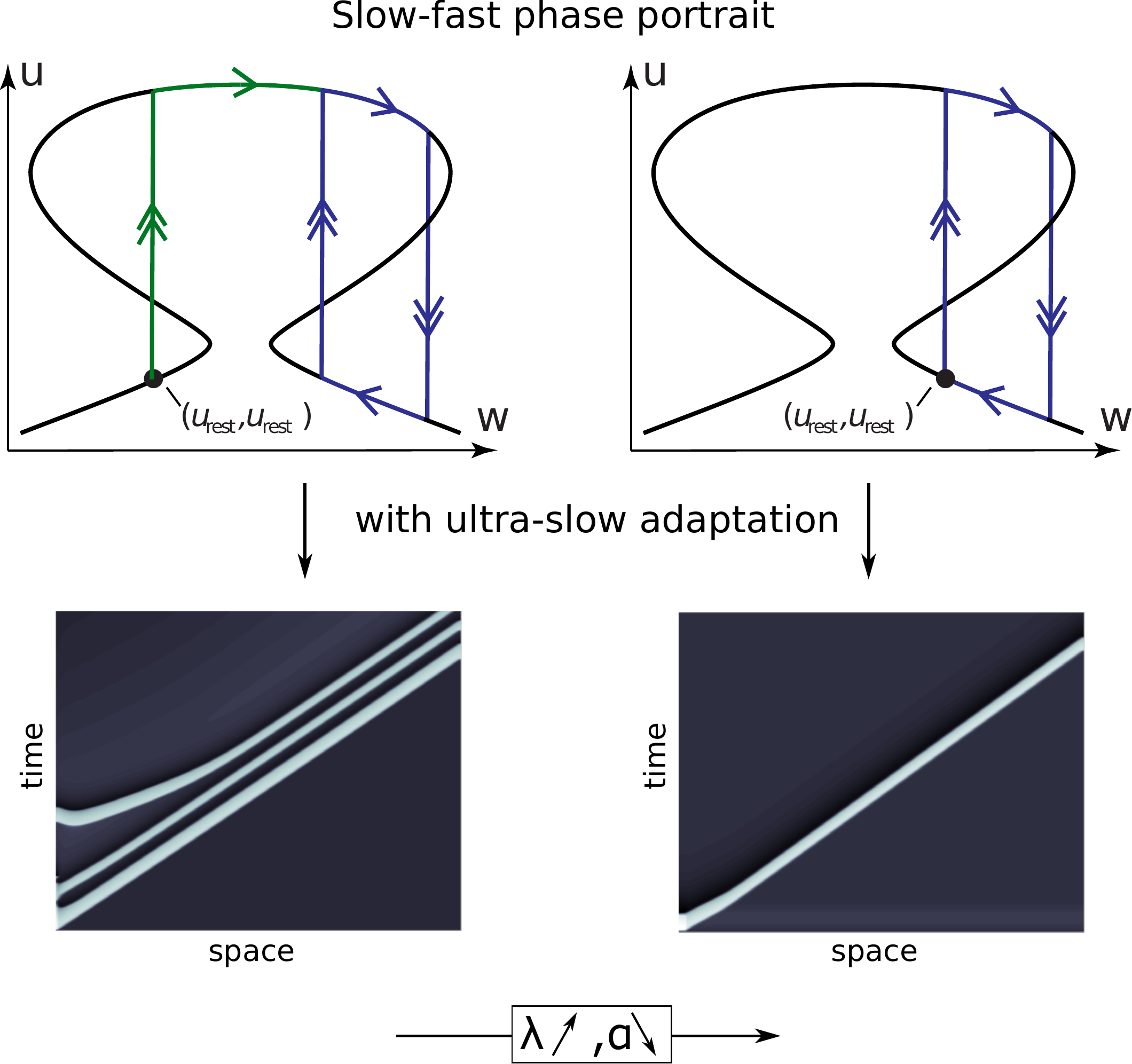}
\caption{Parameters for traveling burst (left plot) are the same as Figure \ref{SimuFig:traveling}. Parameters for traveling pulse (right plot) are the same as traveling bursts except $\lambda=\lambda_{PF}(\beta)+0.15$ and $\alpha=\alpha_{PF}(\beta)+0.4$
}\label{Fig:modulation traveling}
\end{figure}

By changing the value of the bifurcation parameter $\lambda$ and of the unfolding parameter $\alpha$ the traveling burst pattern predicted by Theorem \ref{THM: traveling burst} deforms to a classical traveling pulse (Figure \ref{Fig:modulation traveling}). This is easily understood in terms of the geometry of the singular limit of the slow-fast subsystem. In particular, in the parameter region where there is no bistability in the slow-fast system, we recover the geometry of the Fitzugh-Nagumo model, as the ultra-slow variable barely affects the traveling pulse.

\begin{figure}
\centering
\includegraphics[width=0.6\textwidth]{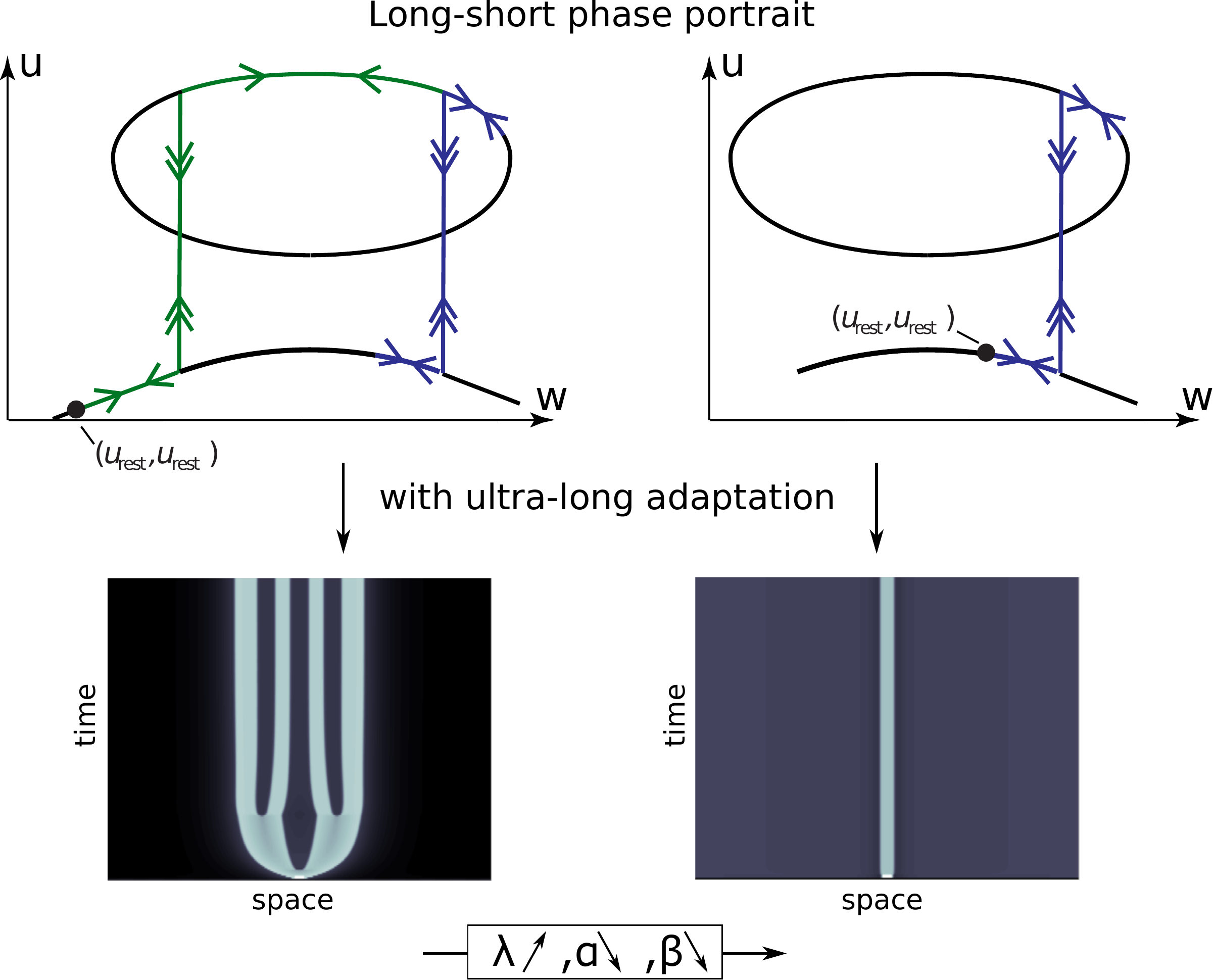}
\caption{Parameters standing burst (left plot) are the same as Figure \ref{SimuFig:standing}. Parameters for standing pulse are the same as standing bursts except $\lambda=\lambda_{PF}(\beta)+0.1$,  $\alpha=\alpha_{PF}(\beta)+0.385$ and $\beta=1/3$.
}\label{Fig:modulation standing}
\end{figure}

The same result holds for the deformation of standing bursts into standing pulses. By changing the value of the bifurcation parameter $\lambda$ and of the unfolding parameters $\alpha$ and $\beta$ the standing burst pattern predicted by Conjecture~\ref{THM: standing burst} is changed into a classical standing pulse (Figure~\ref{Fig:modulation standing}). Again, inspection of the geometry of the singular limit of the short-long range subsystem reveals the qualitative analogy with the geometry of the classical standing pulse solution and the presence of the ultra-long-range variable barely affects the standing pulse solution.

%Changing parameters in different ways it is possible to tune qualitative and quantitative properties of spatio-temporal patterns in model (\ref{EQWcRd3D_temp}) in a finer way. Inspection of the geometry of the singular limit makes this tuning constructive.

We stress that the traveling and standing bursts solutions were found by enforcing the conditions of Theorem~\ref{THM: traveling burst}~and Conjecture~\ref{THM: standing burst}, respectively. The traveling and standing pulse solutions were found by changing the model singular geometry (via modulation of bifurcation and unfolding parameters) to enforce the same qualitative geometry as \cite[Sections 4.2, 4.5, and 5.3]{Jones1995a} and \cite[Sections 4 and 5]{Jones1998},\cite[Section 2]{Dockery1992}, respectively. It would be of interest to perform a numerical continuation analysis to analyze the continuous deformation from one solution to the other.

\section{Discussion}
\label{SEC: Discussion}

\subsection{The behavioral relevance of three-scale bursting models}

In the continuation of the (ode) bursting model introduced in \cite{Franci2014}, 
the (pde) bursting model in the present paper is a three scale model, 
that is, it uses three distinct scales to model the three distinct scales of a bursting wave. 
This complexity may seem unnecessary as bursting has often be described in two
scale reduced ode models (see, e.g., \cite[Chapter 5]{Ermentrout2010a} and references therein). We emphasize, however, that
the robustness and modulation properties discussed here and in \cite{Franci2014} result from
the geometric organization of the attractor, a property that is not shared by two scale 
bursting models.

Two-scale bursting models rely on a 2-dimensional fast subsystem that contains
a family of periodic orbits ending in a homoclinic bifurcation (see, e.g., \cite{Lee1999}) or fold limit cycle \cite{Su2003}. 
The existence of such a family may hold for a range of parameters, but it 
is not a geometric property of the fast system. In contrast, in the three-scale model,
the existence in the slow-fast subsystem of a family of periodic orbits and the associated global
bifurcations proceed from a geometric construction (see \cite[Section 6]{Franci2012} and \cite[Lemma 6]{Franci2014}).
The geometric construction of the global singular skeleton allows for a mapping
between the algebraic structure of the universal unfolding of the organizing singularity and the observed
dynamical behavior. In other words, close to the three-time-scale singular limit, the observed
dynamical behavior is fully determined by the algebraic property of the organizing center. Such a mapping
does not exist in two-time-scale bursting models. 

The mapping between the algebraic structure and the dynamical behavior is key to a 
rigorous study of modulation of three-scale spatio-temporal patterns in terms of paths in the universal
unfolding of the organizing singularity. The physiological relevance of this analysis in the purely
temporal case is illustrated in \cite[Section 3]{Franci2014} and \cite{Drion2015}. We anticipate a similar potential
for the three scale model of the present paper in the analysis of spatio-temporal pattern formation.

\subsection{Stability of bursting waves}
\label{SEC: stab conjecture}

Our work establishes the existence, not the stability, of traveling and standing bursting waves. The stability analysis is beyond the scope of the present paper but predictions can be made in accordance with what is known about the stability of the two-scale models with a cubic nonlinearity. Standing pulses in two-scale models are known to be stable if the space-scale separation $\delta_s$ is much smaller than the time-scale separation \cite[Theorem 4.2]{Rubin1999}, namely
$$\delta_l \ll \varepsilon_s$$
Likewise, to the best of our knowledge, the stability of traveling pulses has been analyzed only in the absence of diffusion in the adaptation variable \cite{Jones1984}, which corresponds to the limit
$$\delta_l =\infty,$$ 
but one can expect that traveling pulses are stable when space-scale separation dominates the time-scale separation, i.e. $$\varepsilon_s\ll  \delta_l.$$

We suspect that similar properties extend to the three-scale model analyzed in the present paper : traveling burst should be stable when time-scale separation dominates space-scale separation, namely
$$\delta_l \ll \varepsilon_s,\text{ and }\delta_{ul}\ll\varepsilon_{us},$$
and standing bursts should be stable when space-scale separation dominates time-scale separation
$$\delta_l \gg \varepsilon_s,\text{ and }\delta_{ul}\gg\varepsilon_{us}.$$

\subsubsection{Numerical illustration}
\label{SEC: numerics}

The following numerical examples were obtained with MATLAB function {\tt pdepe} with no-flux boundary conditions. Details about the used meshes and parameters are given in the figure captions.

{\bf Traveling burst}

\begin{figure}
\centering
\includegraphics[width=0.7\textwidth]{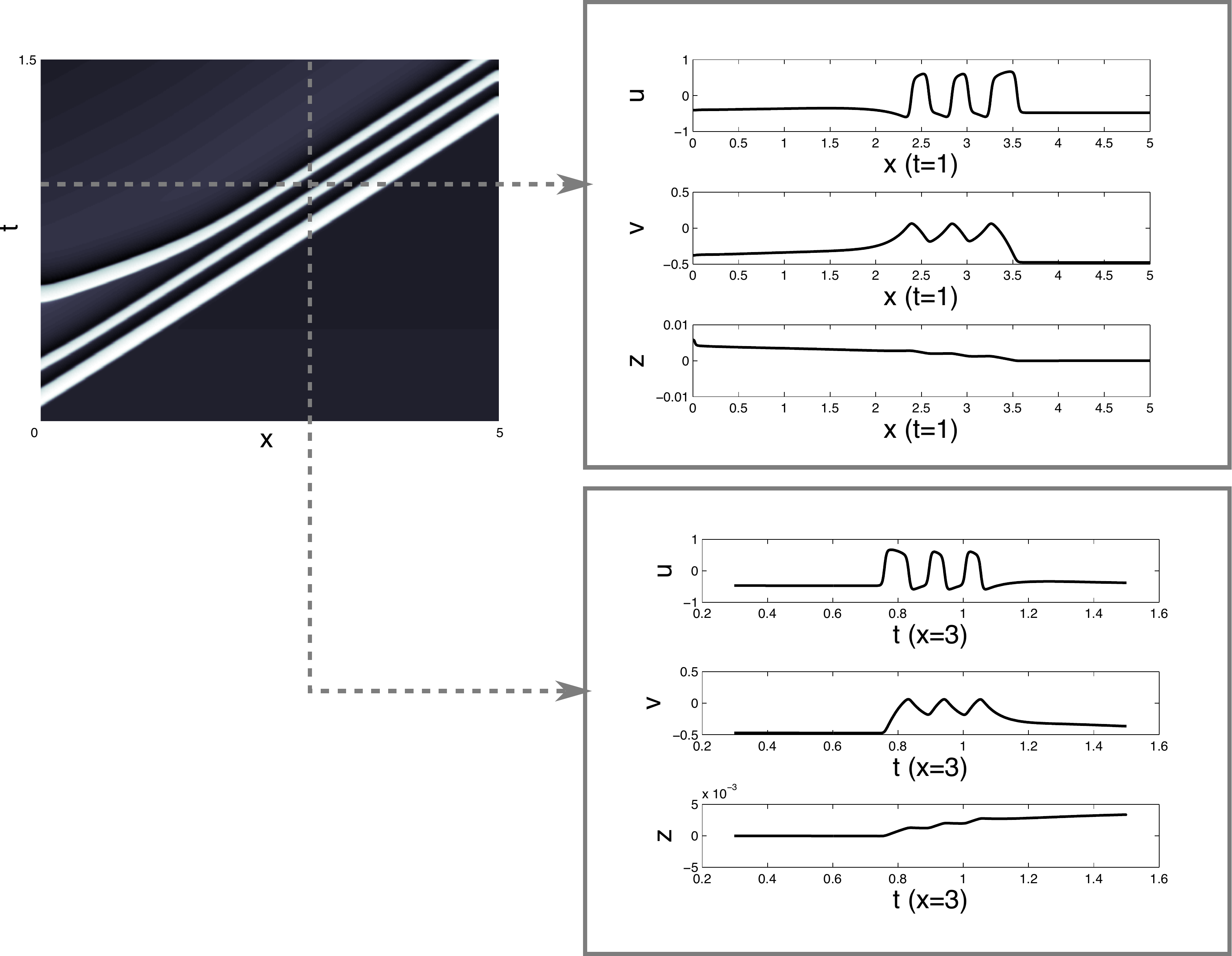}
\caption{Space mesh: $[0:5/3000:5]$, time mesh: $[0:5/1000:5]$.  Time constants: $\tau_u=0.001,\ \tau_w=0.1,\ \tau_z=60$. Diffusion constants: $D_u=0.00005,\ D_w=0.0,\ D_z=0.0$. Parameters: $\beta=1/3$, $\gamma=\gamma_{PF}(\beta)$,  $\lambda=\lambda_{PF}(\beta)-0.02$, $\alpha=\alpha_{PF}(\beta)+0.47$, where the functions $\lambda_{PF}(\cdot),\alpha_{PF}(\cdot),\gamma_{PF}(\cdot)$ are defined in \cite[Appendix A]{Franci2014}. The traveling burst was elicited by perturbing the homogeneous resting state with a perturbation $Pert(t,x)=I_{\{0<x<0.05\}}\, I_{\{0<t<0.25\}}$, where $I_{\{\cdot\}}$ is the indicator function.}\label{SimuFig:traveling}
\end{figure}

Conditions of Theorem \ref{THM:Wavefront} provides a constructive way to chose bifurcation and unfolding parameters as well as time constants in (\ref{EQWcRd3D_temp}) to observe traveling bursting waves. Figure \ref{SimuFig:traveling} shows the result of the numerical integration for $D_w=D_z=0$.

\begin{figure}
\centering
\includegraphics[width=0.7\textwidth]{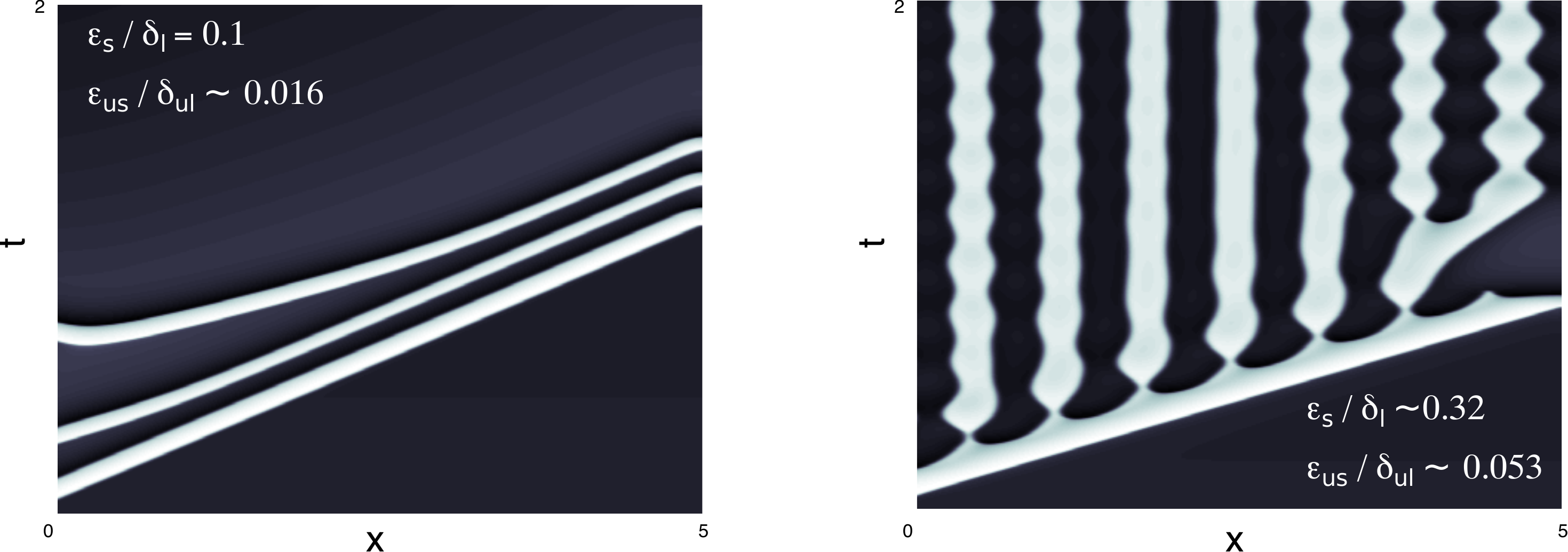}
\caption{Same parameters as Figure \ref{SimuFig:traveling} except $D_w=0.01,D_z=1.0$ in the left plot and  $D_w=0.1,D_z=100.0$ in the right plot.
}\label{SimuFig:traveling robu}
\end{figure}

Figure \ref{SimuFig:traveling robu} shows simulation for non-zero $D_w$ and $D_z$. The simulations support the prediction of Section \ref{SEC: stab conjecture}, that is, it is the ratio between time-scale separation and space-scale separation that determines the stability of the traveling bursting wave. When this ratio is small, that is, time-scale separation dominates space-scale separation, the traveling burst is stable. As this ratio increases and space-scale separation becomes more relevant, the traveling burst loses stability. The resulting pattern exhibits a  ``breathing'' behavior, that is, standing pulses of variable width.

{\bf Standing bursts}

\begin{figure}
\centering
\includegraphics[width=0.7\textwidth]{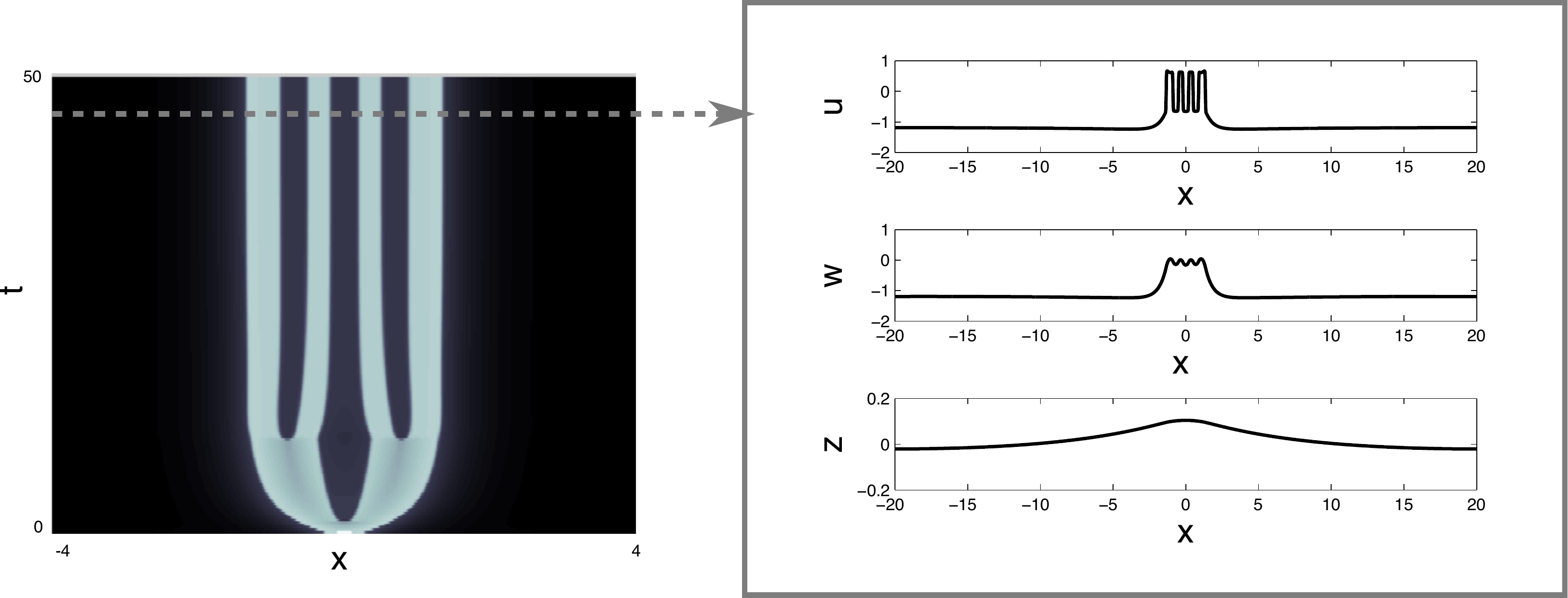}
\caption{Space mesh: $[-20:20/6000:20]$, time mesh: $[0:50/1000:50]$.  Time constants: $\tau_u=\tau_w=\tau_z=0.01$. Diffusion constants: $D_u=0.0001,\ D_w=0.025,\ D_z=100$. Parameters: $\beta=1/3+0.1$, $\gamma=\gamma_{PF}(\beta)$,  $\lambda=\lambda_{PF}(\beta)-0.275$, $\alpha=\alpha_{PF}(\beta)+1.125$ The traveling burst was elicited by perturbing the homogeneous resting state with a perturbation $Pert(t,x)=I_{\{-0.1<x<0.1\}}\, I_{\{0<t<0.5\}}$}\label{SimuFig:standing}
\end{figure}

Conditions of Conjecture \ref{THM: standing burst} provide a constructive way to chose bifurcation and unfolding parameters as well as time constants in (\ref{EQWcRd3D_temp}) to observe traveling bursting waves. Figure \ref{SimuFig:standing} shows the result of the numerical integration in the absence of time-scale separation. In line with the singular construction of Section \ref{SSEC: standing bursts}, we observe a switch from the singular limit (\ref{EQWcRd3DStandODElayerA}),(\ref{EQWcRd3DStandODEreducedA}) to the singular limit (\ref{EQWcRd3DStandODElayerB}),(\ref{EQWcRd3DStandODEreducedB}) at the transition from the quasi-steady part of the pattern to the oscillatory part: $w\sim u$ in the quasi-steady part of the standing burst pattern, whereas only $u$ is at quasi steady state in the oscillatory part.

\begin{figure}
\centering
\includegraphics[width=0.85\textwidth]{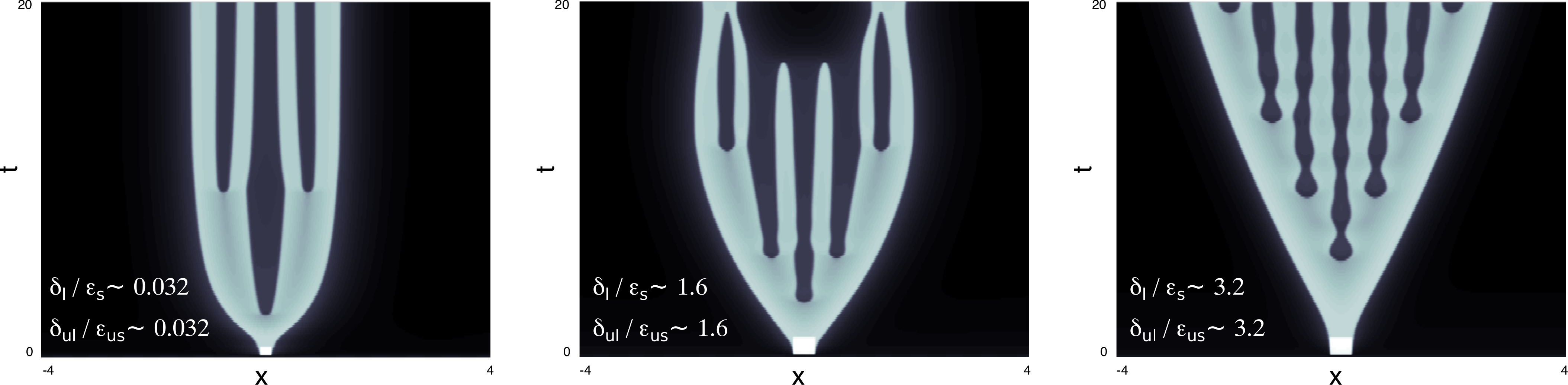}
\caption{Same parameters as Figure \ref{SimuFig:standing} except $\tau_w=0.1,\tau_z=1.0$ in the left plot,  $\tau_w=0.5,\tau_z=25.0$ in the center plot, and $\tau_w=1.0,\tau_z=100.0$ in the right plot.
}\label{SimuFig:standing robu}
\end{figure}

Figure \ref{SimuFig:standing robu} shows the simulation for increasing time-scale separation. Here again, the numerical simulations support the prediction that it is the dominance of the space-scale separation over the time-scale separation that determines the stability of the standing pulse. For small time-scale separations, the pattern exhibits the same standing burst. Increasing time-scale separation the standing burst loses stability. For even larger time-scale separation the pattern begins to exhibit the same breathing behavior observed in Figure \ref{SimuFig:traveling robu}.

\subsection{Multi-scale geometric and stability analysis}

The three-scale model proposed in this paper might motivate an extension of Fenichel theory to more than two scales. To the best of our knowledge, the only general effort in this direction can be found in \cite{Cardin2014} (but see also \cite[Section 6]{Vo2013} and reference therein and \cite{Nan2015} for less rigorous treatments). The reader will indeed notice that the construction of the three-scale standing burst requires two distinct singular limits and, in both limits, only normally hyperbolic part of the critical manifold are visited. This is in contrast to other multi-scale phenomena that are still amenable to the classical two-scale theory of Fenichel (eventually in conjunction with desingularization techniques), such as mixed mode oscillations \cite{Desroches2012,Krupa2008,Nan2015} and non-classical relaxation oscillations \cite{Kuehn2015}. Beside desingularization, the analysis in those reference is in spirit similar to the traveling burst analysis of the present paper, in particular, in terms of conditions imposed on the timescale ratios that allow to reduce the analysis to a two-scale singularly perturbed dynamics. The extension of the Exchange Lemma to general multi-scale singularly perturbed dynamical systems seems of particular relevance.

Also the problem of stability of the constructed multi-scale spatio-temporal pattern should be addressed carefully and more rigorously. The Evan's function technique was successfully applied in the two-scale scenario both for traveling \cite{Jones1984} and standing \cite{Jones1998} pulses. An extension to the three-scale scenario might be natural.

\appendix
\numberwithin{equation}{section}
%\numberwithin{lemma}{section}

\section{Proof of Theorem \ref{THM:Wavefront}}
\label{SSEC: trav front proof}

\begin{figure}[h!]
\center
\includegraphics[width=0.9\textwidth]{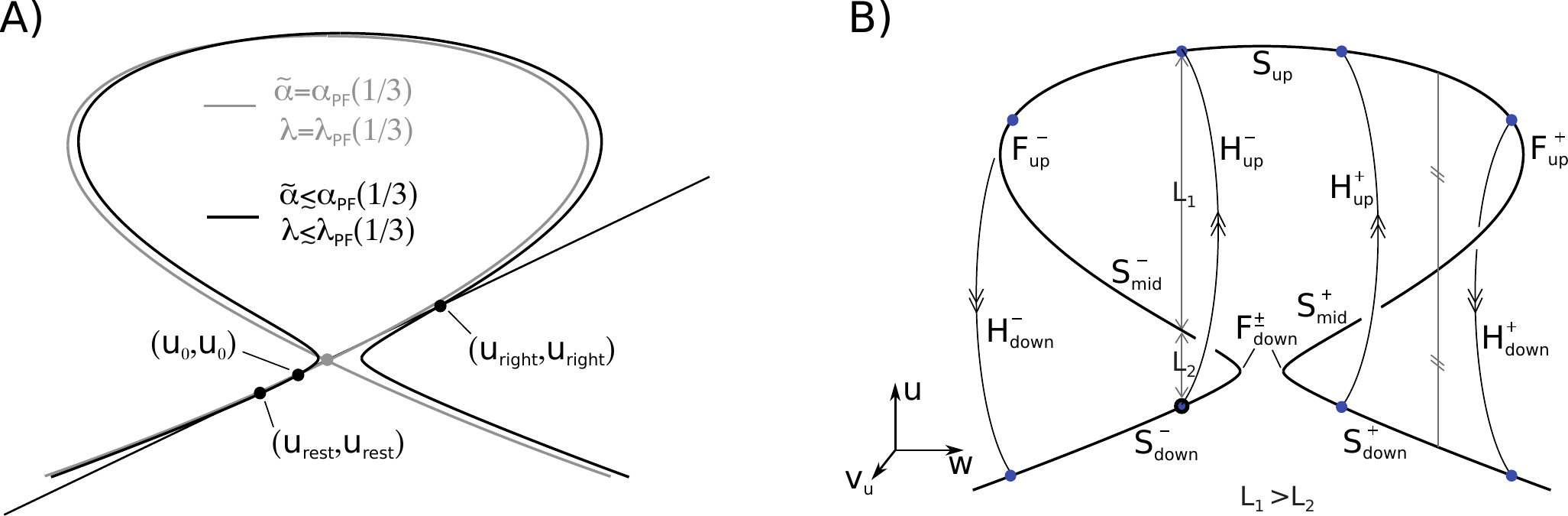}
\caption{{\bf Geometric construction of the singular phase portrait in Figure \ref{FIG3}-A.} Left: critical manifold and fixed points at the pitchfork singularity (gray) and for $\lambda$ and $\alpha$ satisfying the conditions of Lemma \ref{LEM: sing hetero connection} (black). Right: heteroclinic orbits of the layer dynamics and the different invariant manifolds involved in their construction. $\mathcal F^{+/-}_{down}$ and $\mathcal F^{+/-}_{up}$ denotes the four fold singularities in the mirrored hysteresis persistent bifurcation diagram. $S_{down}^{+/-}$, $S_{mid}^{+/-}$, and $S_{up}$ are the disconnected open submanifold of $S$, such that $\mathcal F^{+/-}_{down}\cup\mathcal F^{+/-}_{up}\cup S_{mid}^{+/-}\cup S_{up}=S$. $L_1$ and $L_2$ denotes the distances $d_H(S_{up}\cap \{w=u_{rest}\},S_{mid}^{-}\cap \{w=u_{rest}\})$ and $d_H(S_{down}^{-}\cap \{w=u_{rest}\},S_{mid}^{-}\cap \{w=u_{rest}\})$, respectively. 
}\label{sFIG: front constr}
\end{figure}

We start by a technical lemma that builds a suitable singular skeleton on which we can apply the Exchange Lemma. We refer to Figure~\ref{sFIG: front constr} for the notation.

\vspace{2.5mm}
\begin{lemma}\label{LEM: sing hetero connection}
Let $\beta=\frac{1}{3}$ and $\gamma=\gamma_{PF}\left(\frac{1}{3}\right)=0$. Let $U$ be a neighborhood of $\left(\lambda_{PF}\left(\frac{1}{3}\right),\alpha_{PF}\left(\frac{1}{3}\right)\right)$ in $\mathbb R^2$. There exists an open set $V\subset U\cap\left\{\lambda<\lambda_{PF}\left(\frac{1}{3}\right),\tilde\alpha<\alpha_{PF}\left(\frac{1}{3}\right)\right\} $ such that, for all $(\lambda,\tilde\alpha)\in V$, the following hold
\begin{enumerate}[label=\alph*)]
\item For all $c\in\mathbb R$, the critical manifold $S$ of (\ref{EQWcRd2DTravODE}) has a mirrored hysteresis shape (strong equivalence class {\sf 2} of the persistent bifurcation diagrams of the winged cusp listed at \cite[page 208]{Golubitsky1985}).
\item Model (\ref{EQWcRd2DTravODE}) has exactly three curves of fixed points, corresponding to three roots $u_{rest}<u_{o}<u_{right}$ of (\ref{EQHomSsEq}). The leftmost curve $\{(u_{rest},0,u_{rest},c),\ c\in\mathbb R\}\subset S^-_{down}\subset S$. The middle one $\{(u_{o},0,u_{o},c),\ c\in\mathbb R\}\subset S^-_{down}\cup S^-_{mid}\cup \mathcal F_{down}^- \subset S$. The rightmost one $\{(u_{right},0,u_{right},c),\ c\in\mathbb R\}\subset S^+_{mid}\subset S$.
Moreover, $u_{rest}$ satisfies
$$u_{rest}^3-\beta u_{rest}-\alpha>-u_{rest}^3+\beta u_{rest}-\alpha>0, $$
and
$$u_{rest}^3-(\beta+1)u_{rest}+\alpha>0.$$
\item Let $d_H(\cdot,\cdot)$ denote the Hausdorff distance. Then $$d_H(S_{down}^{-}\cap \{w=u_{rest}\},S_{mid}^{-}\cap \{w=u_{rest}\})<d_H(S_{up}\cap \{w=u_{rest}\},S_{mid}^{-}\cap \{w=u_{rest}\}).$$
\item There exists $c^*\in\mathbb R$ such that the layer dynamics
\begin{IEEEeqnarray}{rCl}\label{sEQ: FHN wcusp ODE layer}
\IEEEyesnumber
u'&=&v\,,\IEEEyessubnumber\\
v'&=&cv-g_{\wc}\left(u,\lambda+w,\tilde\alpha,\frac{1}{3},0\right)\,,\IEEEyessubnumber\\
w'&=&0\,,\IEEEyessubnumber\\
c'&=&0\,,\IEEEyessubnumber
\end{IEEEeqnarray}
of (\ref{EQWcRd2DTravODE}) restricted to the hypersurface $\{c=c^*\}$ has four heteroclinic orbits:
\begin{itemize}
\item $\mathcal H_{up}^{-}$, with base point on the equilibrium $(u_{rest},0,u_{rest})$ and landing point on $S_{up}$;
\item $\mathcal H_{down}^{-}$, with base point on $S_{up}$ at $(u,v,w)=(-u_{rest},0,-\lambda-\sqrt{u_{rest}^3-\beta u_{rest}-\tilde\alpha})$ and landing point on $S_{down}$;
\item $\mathcal H_{up}^{+}$ and $\mathcal H_{down}^{+}$, which are obtained by symmetry of the layer dynamics with respect to the hypersurface $\{w=-\lambda\}$.
\end{itemize}
Let $\bar S_{up}$ and $\bar S_{down}^+$ be compact, connected, normally hyperbolic submanifolds of $S_{up}$ and $S_{down}^+$, respectively, that contain all the base and landing points of the heteroclinic orbits. Then $\mathcal H_{up}^{-}$ is obtained as the transverse intersection (in the $(u,v,w,c)$ space) of the 2-dimensional unstable manifold $W^u_{rest}$ of the curve of fixed points $\{(u_{rest},0,u_{rest},c),\ c\ {\rm near}\ c^*\}$ with the 3-dimensional stable manifold $W^s\left(\bar S_{up}\right)$. The heteroclinic orbits $\mathcal H_{down}^{-},\ \mathcal H_{up}^{+},\ \mathcal H_{down}^{+}$ are obtained as the transverse intersection (in the $(u,v,w,c)$ space) of the 2-dimensional unstable manifold $W^u\left(\bar S_{base}|_{c=c^*}\right)$ of the invariant manifold $\bar S_{base}$ containing the base points restricted to the hypersurface $\{c=c^*\}$ with the 3-dimensional stable manifold $W^s\left(\bar S_{land}\right)$ of the invariant manifold $\bar S_{land}$ containing the landing points.
\end{enumerate}
\end{lemma}
\vspace{2.5mm}

{\bf Proof of Lemma \ref{LEM: sing hetero connection}.}
Points {\it a)}, {\it b)}, and {\it c)} follow from phase plane analysis (Figure \ref{sFIG: front constr}) quantitatively supported by the inspection of transition varieties and persistent bifurcation diagrams of the critical manifold ( $g_{\wc}(u,\lambda,\tilde\alpha,\beta,\gamma)=0$ ) and fixed point (\ref{EQHomSsEq}) equations of (\ref{EQWcRd2DTravODE}), both cubic universal unfolding of the winged cusp, near the singularity at $(u,\lambda,\tilde\alpha,\beta,\gamma)=(-\frac{1}{3},\frac{1}{3},-\frac{2}{27},\frac{1}{3},0)$. This singularity is transcritical for the critical manifold equation and pitchfork for the fixed point equation. For the algebraic expressions of the transition and bifurcation varieties, see \cite[Page 206]{Golubitsky1985} for the critical manifold equation and \cite{Franci2014} for the fixed point equation.

%To prove {\it d)}, we note that, in a neighborhood of the hypersurface $\{w=u_{rest}\}$, (\ref{sEQ: FHN wcusp ODE layer}) is equivalent to the  layer dynamics of the FitzHugh-Nagumo traveling pulse equation \cite[Equation (4.2)]{Jones1995a} with, by point {\it c)}, $0<a<\frac{1}{2}$. By following \cite[Section 4.5]{Jones1995a}, we can therefore conclude the existence of $c^*\neq 0$ such that (\ref{sEQ: FHN wcusp ODE layer}) possesses, for $c=c^*$, the heteroclinic orbit $\mathcal H_{up}^{-}$ satisfying the transversality conditions of point {\it d)}. Existence of the heteroclinic orbits $\mathcal H_{down}^{-}$, $\mathcal H_{up}^{+}$, and $\mathcal H_{down}^{+}$, and their transversality properties follows again by invoking local equivalence of (\ref{sEQ: FHN wcusp ODE layer}) with the layer dynamics of the FitzHugh-Nagumo traveling pulse equation \cite[Equation (4.2)]{Jones1995a} and by following the same construction as in \cite[Section 5.3]{Jones1995a}.

To prove {\it d)} we use existing results on the FitzHugh-Nagumo traveling pulse equation \cite[Sections 4.2 and 5.3]{Jones1995a}. Points {\it a-c)} imply that there exists a diagonal diffeomorphism from a neighborhood of the hypersurface $\{w=u_{rest}\}$ to a neighborhood of the hypersurface $\{w=0\}$, which is the identity on $v$ and is affine in $u$, and which maps the layer dynamics (\ref{sEQ: FHN wcusp ODE layer}) to the layer dynamics of the FitzHugh-Nagumo traveling pulse equation \cite[Eq. (4.2)]{Jones1995a} with parameter $0<a<1/2$ given by
$$a=\frac{d_H(S_{down}^{-}\cap \{w=u_{rest}\},S_{mid}^{-}\cap \{w=u_{rest}\})}{d_H(S_{up}\cap \{w=u_{rest}\},S_{down}^{-}\cap \{w=u_{rest}\})}.$$
Explicitly, the diffeomorphism is given by
\begin{equation}\label{sEQ: diffeo to FHN}
\left(\begin{array}{c}
u\\v\\w\\c
\end{array}\right)\mapsto
\left(\begin{array}{c}
Cu+u_{rest}\\
v\\
-\sqrt{w-u_{rest}^3+\beta u_{rest}-\alpha}-\lambda\\
c
\end{array}\right)\,,
\end{equation} 
where $C$ is a scaling factor such that $u=1$ is the largest of the three roots of transformed layer dynamics fixed point equation computed at $w=0$, that is, $-(Cu+u_{rest})^3+\frac{1}{3}(Cu+u_{rest})-(u_{rest}+\lambda)^2-\alpha=0$. The smallest root is by construction at $u=0$ and, by {\it c)}, the middle root is $u=a$ with $0<a<1/2$. Because by {\it b)} $-u_{rest}^3+\beta u_{rest}-\alpha>0$, the diffeomorphism is well defined.

We now invoke the fact that (\ref{sEQ: diffeo to FHN}) is diagonal, its affine dependence on $u$, and the fact that it is the identity on $v$ and $c$. By the chain rule, these properties implies that the differential forms that we use to track invariant manifolds of (\ref{sEQ: FHN wcusp ODE layer}) and their intersections are given by linear scaling of the same computations as in the classical FitzHugh-Nagumo equation \cite[Section 4.5]{Jones1995a}. As for the classical FitzHugh-Nagumo equation, we conclude the existence of $c^*\neq 0$ such that (\ref{sEQ: FHN wcusp ODE layer}) possesses for $c=c^*$ the heteroclinic orbit $\mathcal H_{up}^{-}$ satisfying the transversality conditions of point {\it d)}.

Existence of the heteroclinic orbits $\mathcal H_{down}^{-}$, $\mathcal H_{up}^{+}$, and $\mathcal H_{down}^{+}$, and their transversality properties follow again from local equivalence of (\ref{sEQ: FHN wcusp ODE layer}) with the layer dynamics of the FitzHugh-Nagumo traveling pulse equation \cite[Equation (4.2)]{Jones1995a} and by following the same construction as in \cite[Section 5.3]{Jones1995a}.

%The existence of the rest of the heteroclinic orbits follows with the same argument but different diffeomorphisms essentially obtained by symmetry arguments. For $\mathcal H_{down}^{-}$, we use antisymmetry of the layer dynamics in $u$ and get
%\begin{equation*}
%\left(\begin{array}{c}
%u\\v\\w\\c
%\end{array}\right)\mapsto
%\left(\begin{array}{c}
%C_1u-u_{rest}\\
%v\\
%-\sqrt{w+u_{rest}^3-\beta u_{rest}-\alpha}-\lambda\\
%c
%\end{array}\right)\,,
%\end{equation*}
%which maps a neighborhood of the hypersurface $\{w=-\sqrt{u_{rest}^3-\beta u_{rest}-\alpha}-\lambda\}$ to a neighborhood of the hypersurface $\{w=0\}$ and the transformed dynamics is again  given by \cite[Eq. (4.2)]{Jones1995a}. Because by {\it b)} $u_{rest}^3-\beta u_{rest}-\alpha>0$, the diffeomorphism is well defined. For $\mathcal H^+_{up/down}$, we use the mirror symmetry of the layer dynamics with respect to the hypersurface $\{w=-\lambda\}$.

\hfill$\square$\\

To prove the theorem, we track the 2-dimensional unstable manifold $W^u_{rest}$ of the curve of fixed points $\{(u_{rest},0,u_{rest},c),\ c\ {\rm near}\ c^*\}$, that is, we follow the mapping of a germ of $W^u_{rest}$ through the flow associated to (\ref{EQWcRd2DTravODE}). For $\varepsilon_s>0$ and sufficiently small, let  $\bar S_{up,\varepsilon_s}$ and $\bar S_{down,\varepsilon_s}^+$ be the slow manifolds obtained as Fenichel perturbation of $\bar S_{up}$ and $\bar S_{down}^+$, respectively. Thanks to Lemma \ref{LEM: sing hetero connection}{\it d)}, we can apply \cite[Lemma 4.1]{Jones1994} and the Exchange Lemma and follow the same arguments as the last two paragraphs of the proof of the Theorem in \cite[Section 4]{Jones1994} to conclude the following: $W^u_{rest}$ intersects transversely $W^s\left(\bar S_{up,\varepsilon_s}\right)$ along $\mathcal H^-_{up}$; between $\mathcal H^-_{up}$ and $\mathcal H^+_{down}$ it lies $\mathcal O(\varepsilon_s)$-close to $S_{up}$. It leaves $S_{up}$ along $\mathcal H^+_{down}$ intersecting transversely $W^s\left(\bar S_{down,\varepsilon_s}^+\right)$. Between $\mathcal H^+_{down}$ and $\mathcal H^+_{up}$ it lies $\mathcal O(\varepsilon_s)$-close to $S_{down}^+$; finally, it leaves $S_{down}^+$ along $\mathcal H^+_{up}$, again, intersecting transversely $W^s\left(\bar S_{up,\varepsilon_s}\right)$. We can continue to track the forward mapping of $W^u_{rest}$ through the flow associate to (\ref{EQWcRd2DTravODE}) to conclude that it transversely intersects $W^s\left(\bar S_{up,\varepsilon_s}\right)$ and $W^s\left(\bar S_{down,\varepsilon_s}^+\right)$ infinitely many times. Since $c$ is a parameter, all the (one-dimensional) transverse intersections are in the same $\{c=\bar c\}$-slice, for some $\bar c=c^*+\mathcal O(\varepsilon_s)$. We claim that the trajectory containing all the transverse intersections, call it $h_{\varepsilon_s}$, is the heteroclinic trajectory corresponding to the traveling front.

It remains to show that, in the same $\{c=\bar c\}$-slice there also exists a periodic solution $\ell_{\varepsilon_s}$ of (\ref{EQWcRd2DTravODE}) and that $h_{\varepsilon_s}\subset W^s(\ell_{\varepsilon_s})$. To this aim, we track the two-dimensional unstable manifold $W^u\left(\bar S_{up,\varepsilon_s}|_{c=\bar c}\right)$. By the Exchange Lemma, $W^u\left(\bar S_{up,\varepsilon_s}|_{c=\bar c}\right)$ and $W^u_{rest}$ are $O(\varepsilon_s)$-$C^1$ close to each other at the heteroclinic jump $\mathcal H^+_{down}$. It follows that, for $\varepsilon_s$ sufficiently small,
%after an $O(1/\varepsilon)$-time near $S_{down}^+$,
also $W^u\left(\bar S_{up,\varepsilon_s}|_{c=\bar c}\right)$ comes back (under the flow associated to (\ref{EQWcRd2DTravODE})) in a neighborhood of $S_{up}$ transversely intersecting $W^s\left(\bar S_{up,\varepsilon}\right)$ inside $\{c=\bar c\}$. 
%By Lemma \ref{LEM: sing hetero connection}{\it c)} there are no fixed points on $\bar S_{up,\varepsilon}$ nor on $\bar S_{down,\varepsilon}^+$.
It follows that the 1-dimensional transverse (in the $(u,v,w,c)$ space) intersection $W^u\left(\bar S_{up,\varepsilon_s}|_{c=\bar c}\right )\cap_T W^s\left(\bar S_{up,\varepsilon_s}\right)$ defines a periodic orbit $\ell_{\varepsilon_s}$ of (\ref{EQWcRd2DTravODE}) in the slice $c=\bar c$. Because $W^s(\ell_{\varepsilon_s})=W^s\left(\bar S_{up,\varepsilon_s}|_{c=\bar c}\right)$ and $h_{\varepsilon_s}\subset W^s\left(\bar S_{up,\varepsilon_s}|_{c=\bar c}\right)$, the statement is proved for $\beta=1/3$ and $\gamma=0$. For $\beta\neq 1/3,\ \gamma\neq0$, the theorem follows from the persistence of all the transverse intersections used in the construction to arbitrary $C^1$ perturbations.\hfill$\square$

\section{Proof of Theorem \ref{THM:Standfront}}
\label{SSEC: stand front proof}

We use two technical lemmas to build a suitable singular skeleton on which we can apply the Exchange Lemma.

\vspace{2.5mm}
\begin{lemma}\label{LEM:stand front constr 1}
Let $\beta=\frac{1}{3}$ and $\gamma=\gamma_{PF}\left(\frac{1}{3}\right)=0$. For all $\tilde\alpha\in[\alpha_{PF}(1/3),0)$ and all $\lambda\in\mathbb R$, the layer dynamics of (\ref{EQWcRd2DStandODE})
\begin{IEEEeqnarray}{rCl}\label{sEQWcRd2DStandODEfastLem}
\IEEEyesnumber
u'&=&v_u\,,\IEEEyessubnumber\\
v_u'&=&-g_{\wc}(u,\lambda+w,\tilde\alpha,\frac{1}{3},0)\,,\IEEEyessubnumber
\end{IEEEeqnarray}
possess two heteroclinic orbits at $w=w_{h_1}(\tilde\alpha)$ and two heteroclinic orbits at $w=w_{h_2}(\tilde\alpha)$, with $w_{h_1}(\tilde\alpha)=-\lambda-\sqrt{-\tilde\alpha}<-\lambda< -\lambda+\sqrt{-\tilde\alpha} = w_{h_2}(\tilde\alpha)$. For $\tilde\alpha=0$, model (\ref{sEQWcRd2DStandODEfastLem}) possesses two heteroclinic orbits at $w=-\lambda$.
\end{lemma}
\vspace{2.5mm}
{\bf Proof of Lemma \ref{LEM:stand front constr 1}.}
For $w=w_{h_1}(\tilde\alpha)$ or $w=w_{h_2}(\tilde\alpha)$, equation (\ref{sEQWcRd2DStandODEfastLem}) reduces to
\begin{IEEEeqnarray*}{rCl}
u'&=&v_u\,,\\
v_u'&=&-u^3+\frac{u}{3}\,,
\end{IEEEeqnarray*}
which is Hamiltonian, has three fixed point at $v_u=0$ and $u=-\frac{1}{\sqrt{3}},0,\frac{1}{\sqrt{3}}$ and satisfies
$$\int_{-\frac{1}{\sqrt{3}}}^{\frac{1}{\sqrt{3}}} \left( -u^3+\frac{u}{3}\right)du =0 $$
\hfill$\square$

\begin{figure}[h!]
\center
\includegraphics[width=0.9\textwidth]{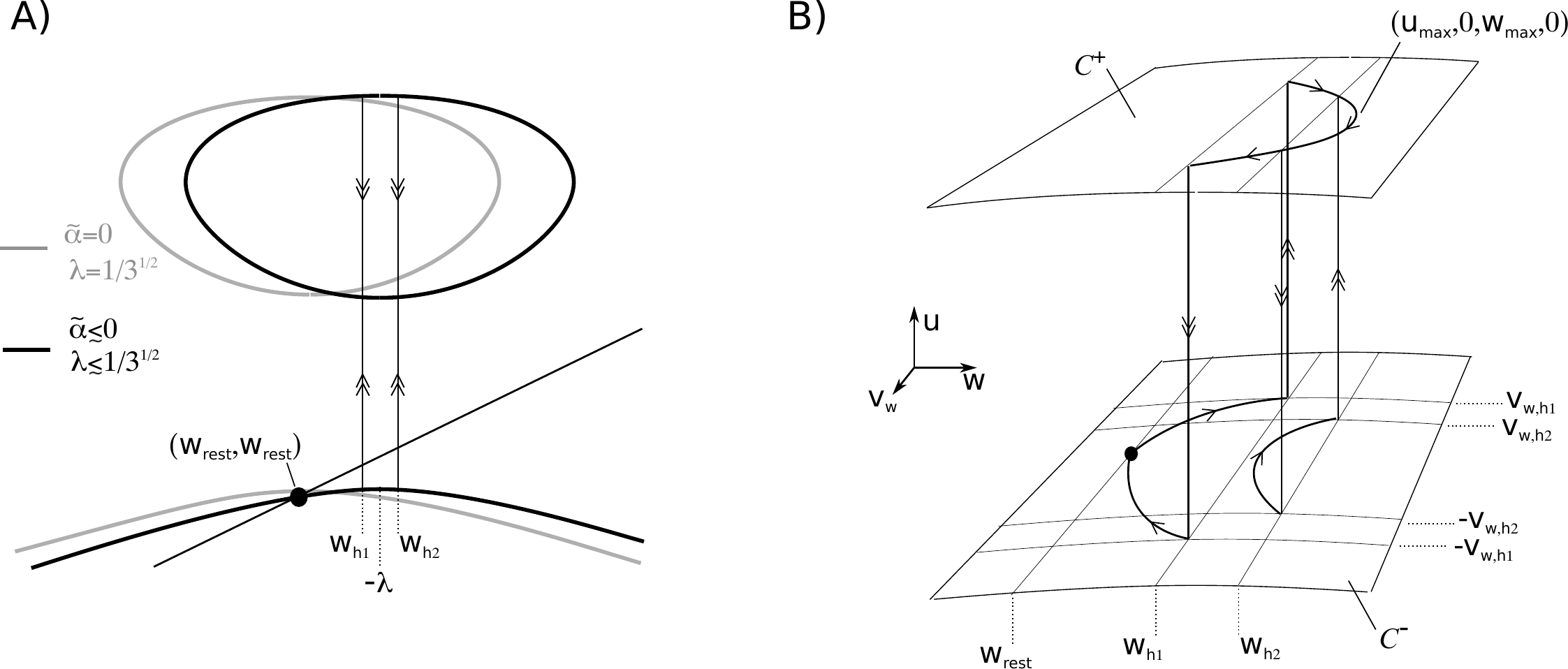}
\caption{{\bf Geometric construction of the singular phase portrait in Figure \ref{FIG3}-B.}}\label{sFIG: stand front constr}
\end{figure}

The notation in the next lemma is defined in Figure \ref{sFIG: stand front constr}.

\vspace{2.5mm}
\begin{lemma}\label{LEM:stand front constr 2}
Let $\beta=\frac{1}{3}$ and $\gamma=\gamma_{PF}\left(\frac{1}{3}\right)=0$. There exist a neighborhood $U$ of $(\lambda,\tilde\alpha)=(1/\sqrt{3},0)$ and an open set $V\subset U\cap\{\lambda<1/\sqrt{3},\,\tilde\alpha<0\}$, such that, for all $(\lambda,\tilde\alpha)\in V$ the following hold:
\begin{enumerate}[label=\alph*)]
\item The critical manifold of model  (\ref{EQWcRd2DStandODE}) belongs to class {\sf 3} of the persistent bifurcation diagrams of the winged cusp \cite[page 208]{Golubitsky1985}. Let $w_{fold}$ be the $w$-coordinate of the right fold of the critical manifold.  Let $u=u_{down}(w)$ be the function such that the lower branch of the critical manifold $\mathcal C^-=\{(u,v_u,w,v_w)=(u_{down}(w),0,w,0)\}$ and $u_{up}(w)$ the function such that the upper branch of the critical manifold between the two fold singularities $\mathcal C^+=\{(u,v_u,w,v_w)=(u_{up}(w),0,w,0)\}$.
\item Model (\ref{EQWcRd2DStandODE}) has a unique fixed point at $(w_{rest},0,w_{rest},0)$, with $w_{rest}<-\lambda$.
%This fixed point is a saddle for the reduced dynamics (\ref{EQWcRd2DStandODEslow}) with a one-dimensional stable and one dimensional unstable manifold.
\item The following inequalities hold:
\begin{IEEEeqnarray}{rCl}\label{sEQ:stand front int cond}
\IEEEyesnumber
\int_{w_{rest}}^{w_{h_1}} \big( u_{down}(s)-s \big) ds &+& \int_{w_{h_1}}^{w_{fold}} \big( u_{up}(s)-s \big) ds >0\IEEEyessubnumber\\
\int_{w_{rest}}^{w_{h_1}} \big( u_{down}(s)-s \big) ds &+& \int_{w_{h_1}}^{w_{h_2}} \big( u_{up}(s)-s \big) ds <0\IEEEyessubnumber
\end{IEEEeqnarray}
%\item Let $W^{u}_{sl}$ be embedding of the unstable manifold of $(w_{rest},0,w_{rest},0)$ in $\mathbb R^4$ and let $W^{cu}$ be the union of the 1-dimensional fast unstable manifolds of points in the 1-dimensional slow unstable manifold $W^{u}_{sl}$. Then 
\end{enumerate}
\end{lemma}
\vspace{2.5mm}
{\bf Proof of Lemma \ref{LEM:stand front constr 2}.}
Points {\it a)}, {\it b)} follow from phase plane analysis (Figure \ref{sFIG: stand front constr}) quantitatively supported by the inspection of transition varieties and persistent bifurcation diagrams of the critical manifold ( $g_{\wc}(u,\lambda,\tilde\alpha,\beta,\gamma)=0$ ) and fixed point (\ref{EQHomSsEq}) equations of (\ref{EQWcRd2DTravODE}), both cubic universal unfolding of the winged cusp, near the singularity at $(u,\lambda,\tilde\alpha,\beta,\gamma)=(x_{PF}(\frac{1}{3}),\lambda_{PF}(\frac{1}{3}),\alpha_{PF}(\frac{1}{3}),1/3,\gamma_{PF}(\frac{1}{3}))=(-\frac{1}{3},\frac{1}{3},-\frac{2}{27},\frac{1}{3},0)$. This singularity is transcritical for the critical manifold equation and pitchfork for the fixed point equation. For the algebraic expressions of the transition and bifurcation varieties, see \cite[Page 206]{Golubitsky1985} for the critical manifold equation and \cite{Franci2014} for the fixed point equation.

Point {\it c)}. Because $\int_{w_{h_1}}^{w_{fold}} \big( u_{up}(s)-s \big) ds >0$ and $w_{rest}=w_{h_1}$ for $(\lambda,\tilde\alpha)=(1/\sqrt{3},0)$, it follows that (\ref{sEQ:stand front int cond}a) is satisfied for $(\lambda,\tilde\alpha)=(1/\sqrt{3},0)$. By continuity of the integral operator, the same holds true for $(\lambda,\tilde\alpha)$ sufficiently close to $(1/\sqrt{3},0)$.

To prove the second inequality, observe that $\int_{w_{rest}}^{w_{h_1}} \big( u_{down}(s)-s \big) ds<0$ and $w_{h_1}=w_{h_2}$, for $\lambda<1/\sqrt{3}$ and $\tilde\alpha=0$. Therefore  (\ref{sEQ:stand front int cond}b) is satisfied for $\lambda<1/\sqrt{3}$ and $\tilde\alpha=0$. By continuity, the same holds for $(\lambda,\tilde\alpha)$ sufficiently close to $(1/\sqrt{3},0)$.
\hfill$\square$

The existence of the homoclinic orbit now follow exactly as \cite[Sections 4 and 5]{Jones1998}. The two models share indeed the exact same geometry (compare Figure \ref{sFIG: stand front constr} right and \cite[Figure 4]{Jones1998}). The sole difference is the absence of a fixed point on $\mathcal C^+$. This only changes the construction of the slow portion of the singular homoclinic  orbit on $\mathcal C^+$. Construction of this portion here follows the same line as in \cite[pages 238 and 239]{Dockery1992}.

The same construction is used for the periodic orbit. The only difference is that the initial slow portion is now defined as the slow trajectory passing through the point $(w_{h_2},v_{w,h2})$ instead as the unstable manifold of the fixed point. Then, all the results of \cite[Sections 4 and 5]{Jones1998} still hold.

The fact that the homoclinic and singular periodic orbits are $\mathcal O(\delta_l)$-close to each other near the point $(u_{max},0,w_{max},0)$ where $w$ reaches its maximum along the homoclinic trajectory follows from the fact that, near $\mathcal C^+$, the homoclinic and the periodic orbits are perturbation of, and hence $\mathcal O(\delta_l)$-close to, the same invariant manifold; likewise, for their stable and unstable manifolds.

\section{Proof of Theorem \ref{THM: traveling burst}}
\label{SSEC: trav burst proof}

\begin{figure}[h!]
\includegraphics[width=0.9\textwidth]{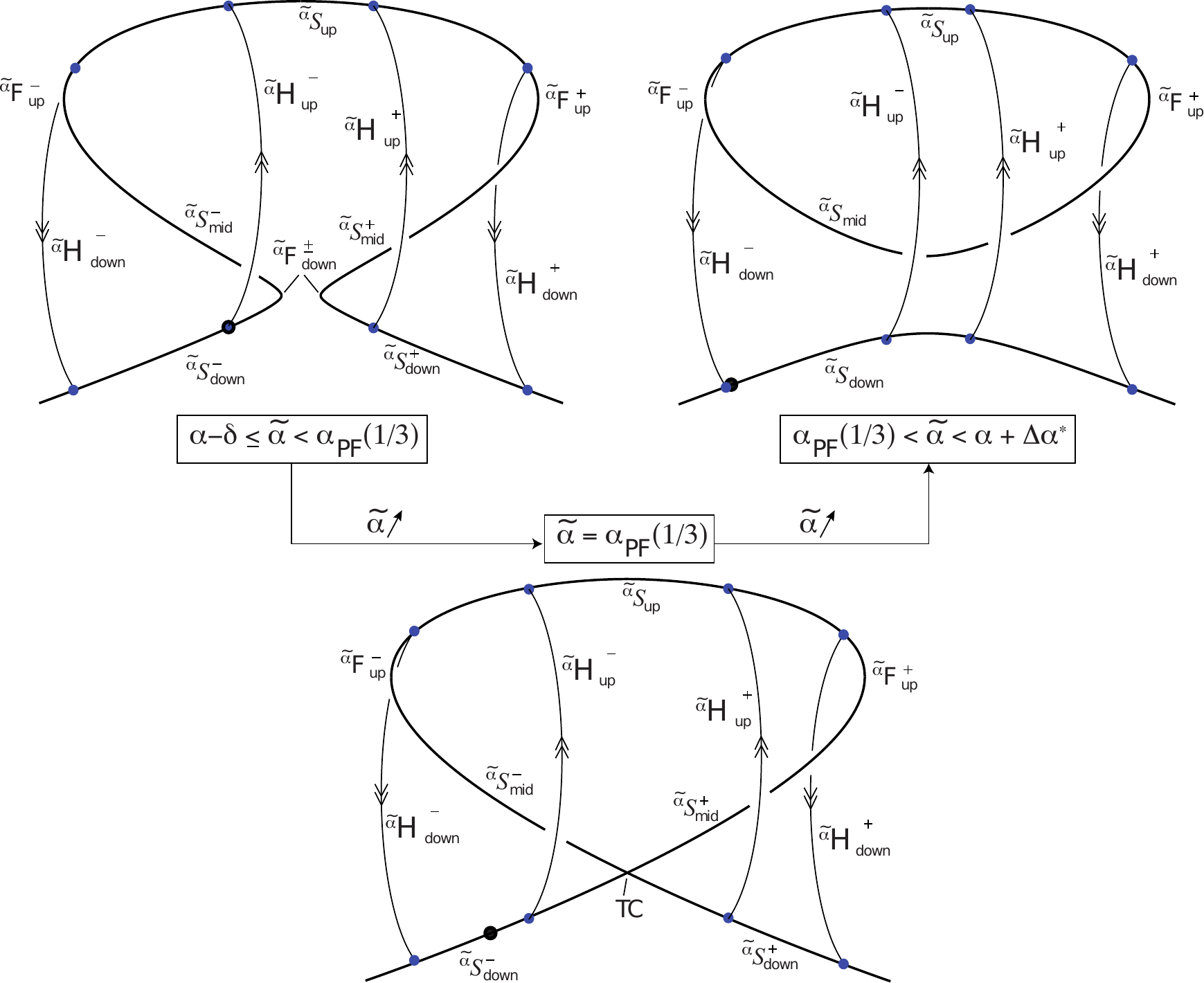}
\caption{}\label{sFIG: burst constr}
\end{figure}

We refer the reader to Figure \ref{sFIG: burst constr} for the notation used in the following lemma.

\vspace{2.5mm}
\begin{lemma}\label{LEM: sing hetero connection 2}
Let $\beta=\frac{1}{3}$ and $\gamma=\gamma_{PF}(\frac{1}{3})=0$. Let $V\subset\mathbb R^2$ and $c^*\neq 0$ be defined as in the statement of Lemma \ref{LEM: sing hetero connection}. For all $(\lambda,\alpha)\in V$, the following hold true
\begin{enumerate}[label=\alph*)]
\item There exists $\Delta\alpha^*>0$ and $\delta>0$ such that, for all $\tilde\alpha\in[\alpha-\delta,\Delta\alpha^*)$ the layer dynamics
\begin{IEEEeqnarray}{rCl}\label{sEQ: FHN wcusp ODE layer2}
\IEEEyesnumber
u'&=&v\,,\IEEEyessubnumber\\
v'&=&cv-g_{\wc}\left(u,\lambda+w,\tilde\alpha,\frac{1}{3},0\right)\,,\IEEEyessubnumber\\
w'&=&0\,,\IEEEyessubnumber\\
c'&=&0\,,\IEEEyessubnumber
\end{IEEEeqnarray}
of (\ref{EQWcRd2DTravODE}) restricted to the hypersurface $\{c=c^*\}$ has four heteroclinic orbits $^{\tilde\alpha}\mathcal H_{up/down}^{+/-}$ that are all obtained as the transverse intersections (in the $(u,v,w,c)$ space) of the 2-dimensional unstable manifold $W^u\left(^{\tilde\alpha}\bar S_{base}|_{c=c^*}\right)$ of the invariant manifold $^{\tilde\alpha}\bar S_{base}$ where the base point lies restricted to the hypersurface $\{c=c^*\}$ with the 3-dimensional stable manifold $W^s\left(^{\tilde\alpha}\bar S_{land}\right)$ of the invariant manifold $^{\tilde\alpha}\bar S_{land}$ where the landing point lies\footnote{As in Lemma \ref{LEM: sing hetero connection}, the overline means a compact, connected, normally hyperbolic submanifold of the relative manifold that contains all the base and landing points of the heteroclinic orbits }.
\item For $\tilde\alpha=\alpha+\Delta\alpha^*$ the two heteroclinic orbits $^{\tilde\alpha}\mathcal H_{up}^{+/-}$ merge in a fast heteroclinic jump $\mathcal H_{up}^0$. This heteroclinic lies in the hypersurface $\{w=-\lambda\}$ and is obtained as the non-transverse intersection $W^u\left(^{\tilde\alpha}\bar S_{down}|_{c=c^*}\right)\cap W^s\left(^{\tilde\alpha}\bar S_{up}\right)$.
\item The two downward heteroclinic jumps $^{\tilde\alpha}\mathcal H_{down}^{+/-}$ persist for $\tilde\alpha\in[\alpha+\Delta\alpha^*,\alpha+\Delta\alpha^*+\delta)$ and with the same transversality properties as in {\it a)}.
\item For all $\tilde\alpha\in[\alpha,\alpha+\Delta\alpha^*)$, let
%and $\mathcal M=\,^{\tilde\alpha}\bar S_{up},\,^{\tilde\alpha}\bar S_{down},\,^{\tilde\alpha}\bar S_{down}^+$, let $\hat{\mathcal M}=\mathcal M\cap\{w\in[^{\tilde\alpha}w_{up},^{\tilde\alpha}w_{down}]\}$, where
$^{\tilde\alpha}w_{up}<^{\tilde\alpha}w_{down}$ be the $w$-coordinate of the heteroclinic $^{\tilde\alpha}\mathcal H_{up}^{+}$ and $^{\tilde\alpha}\mathcal H_{down}^{+}$, respectively. There exists $C>0$ such that, for all $\tilde\alpha\in[\alpha,\alpha_{PF}(1/3)]$,
$$\int_{^{\tilde\alpha}w_{up}}^{^{\tilde\alpha}w_{down}}(u|_{^{\tilde\alpha}S_{down}^+}-(\tilde\alpha-\alpha))dw+\int_{^{\tilde\alpha}w_{up}}^{^{\tilde\alpha}w_{down}}(u|_{^{\tilde\alpha}S_{up}}-(\tilde\alpha-\alpha))dw>C$$
and, for all $\tilde\alpha\in(\alpha_{PF}(1/3),\alpha+\Delta\alpha^*)$,
$$\int_{^{\tilde\alpha}w_{up}}^{^{\tilde\alpha}w_{down}}(u|_{^{\tilde\alpha}S_{down}}-(\tilde\alpha-\alpha))dw+\int_{^{\tilde\alpha}w_{up}}^{^{\tilde\alpha}w_{down}}(u|_{^{\tilde\alpha}S_{up}}-(\tilde\alpha-\alpha))dw>C.$$
\end{enumerate}
\end{lemma}
\vspace{2.5mm}

{\bf Proof of Lemma \ref{LEM: sing hetero connection 2}}.
%Point {\it a)} follows again by local equivalence of the layer dynamics (\ref{sEQ: FHN wcusp ODE layer2}) to the  layer dynamics of the standard FitzHugh-Nagumo pulse equation \cite[Eq. (4.2)]{Jones1995a} and by following the same argument as \cite[Section 5.3]{Jones1995a}.

To prove point {\it a)} we use again local equivalence to the layer dynamics of the standard FitzHugh-Nagumo pulse equation. The diffeomorphism (\ref{sEQ: diffeo to FHN}) between the layer dynamics (\ref{sEQ: FHN wcusp ODE layer2}) and the FitzHugh-Nagumo pulse equation \cite[Eq. (4.2)]{Jones1995a} with the same parameter $a$ as above, can be build also for $\tilde\alpha\neq\alpha$. For $^{\tilde\alpha}\mathcal H_{up}^-$ it reads
\begin{equation*}
\left(\begin{array}{c}
u\\v\\w\\c
\end{array}\right)\mapsto
\left(\begin{array}{c}
C_{\tilde\alpha}u+u_{rest}\\
v\\
-\sqrt{w-u_{rest}^3+\beta u_{rest}-\tilde\alpha}-\lambda\\
c
\end{array}\right)\,,
\end{equation*}
which maps a neighborhood of the hypersurface $\{w=-\sqrt{-u_{rest}^3+\beta u_{rest}-\tilde\alpha}-\lambda\}$ to a neighborhood of $\{w=0\}$, whereas for $^{\tilde\alpha}\mathcal H_{down}^-$ it reads
\begin{equation*}
\left(\begin{array}{c}
u\\v\\w\\c
\end{array}\right)\mapsto
-\left(\begin{array}{c}
C_{\tilde\alpha}u-u_{rest}\\
v\\
\sqrt{w+u_{rest}^3-\beta u_{rest}-\tilde\alpha}-\lambda\\
c
\end{array}\right)\,,
\end{equation*}
which maps a neighborhood of the hypersurface $\{w=-\sqrt{u_{rest}^3-\beta u_{rest}-\tilde\alpha}-\lambda\}$ to a neighborhood of $\{w=0\}$. For $^{\tilde\alpha}\mathcal H_{up/down}^+$ we use symmetry with respect to the hypersurface $\{w=-\lambda\}$. The $w$-coordinate of the heteroclinics $^{\tilde\alpha}\mathcal H_{up}^-$ and $^{\tilde\alpha}\mathcal H_{down}^-$ are given by $^{\tilde\alpha}w_{up}^-=-\sqrt{-u_{rest}^3+\beta u_{rest}-\tilde\alpha}-\lambda$ and $^{\tilde\alpha}w_{down}^-=-\sqrt{u_{rest}^3-\beta u_{rest}-\tilde\alpha}-\lambda$, respectively. Their mirrors with respect to the hypersurface $\{w=-\lambda\}$ are $^{\tilde\alpha}w_{up}^+=\sqrt{-u_{rest}^3+\beta u_{rest}-\tilde\alpha}-\lambda$ and $^{\tilde\alpha}w_{down}^-=\sqrt{u_{rest}^3-\beta u_{rest}-\tilde\alpha}-\lambda$. Existence of the four heteroclinic orbits and their transversality conditions then follows as in \cite[Section 5.3]{Jones1995a}.

The value $\Delta\alpha^*$ of point {\it b)} is computed by imposing $-u_{rest}^3+\beta u_{rest}-\alpha-\Delta\alpha^*=0$. As $\tilde\alpha\to\alpha+\Delta\alpha^*$, the two heteroclinics $^{\tilde\alpha}\mathcal H_{up}^-$ and $^{\tilde\alpha}\mathcal H_{up}^-$ converge to each other in the Hausdorff distance, because $\lim_{\tilde\alpha\to\alpha+\Delta\alpha^*}\, ^{\tilde\alpha}w_{up}^-=\lim_{\tilde\alpha\to\alpha+\Delta\alpha^*}\, ^{\tilde\alpha}w_{up}^+=-\lambda$. For $\tilde\alpha=\alpha+\Delta\alpha^*$, the map
$$\left(\begin{array}{c}
u\\v\\-\lambda\\c
\end{array} \right)\mapsto\left(\begin{array}{c}
Cu+u_{rest}\\
v\\
0\\
c
\end{array}\right) $$
maps the layer dynamics (\ref{sEQ: FHN wcusp ODE layer2}) restricted to the hypersurface $\{w=-\lambda\}$ to the traveling wave problem \cite[Eq. (4.10)]{Jones1995a} with the same parameter $a$ as above. With the same computation as \cite[Section 4.5]{Jones1995a}, we conclude the existence of the heteroclinic $\mathcal H_{up}^0$ in the same hypersurface $\{c=c^*\}$ as $^{\tilde\alpha}\mathcal H^{+/-}_{up/down}$, $\alpha\in[\alpha,\alpha+\Delta\alpha)$. The non-transversality condition follows by the fact that transversality is not compatible with a (local) change in the number of transverse intersections.

By Lemma \ref{LEM: sing hetero connection}{\it b)}, we have $u_{rest}^3-\beta u_{rest}>-u_{rest}^3+\beta u_{rest}$. This in turn implies that the outer heteroclinics $\mathcal H_{down}^{+/-}$ persist and with same properties as in {\it a)} for $\tilde\alpha=\alpha+\Delta\alpha^*$ and, by continuity, also for $\alpha+\Delta\alpha^*<\alpha<\alpha+\Delta\alpha^*+\delta$, which proves {\it c)}.

To prove point {\it d)}, we start by noticing that, by antysimmetry in $u$ of the layer dynamics (\ref{sEQ: FHN wcusp ODE layer}), for all $\tilde\alpha\in[\alpha,\alpha_{PF}(\frac{1}{3})]$,
$$\int_{^{\tilde\alpha}w_{up}}^{^{\tilde\alpha}w_{down}}u|_{^{\tilde\alpha}S_{down}^+} dw = - \int_{^{\tilde\alpha}w_{up}}^{^{\tilde\alpha}w_{down}}u|_{^{\tilde\alpha}S_{up}} dw$$
and, for all $\tilde\alpha\in(\alpha_{PF}(\frac{1}{3}),\alpha+\Delta\alpha^*)$,
$$\int_{^{\tilde\alpha}w_{up}}^{^{\tilde\alpha}w_{down}}u|_{^{\tilde\alpha}S_{down}} dw = - \int_{^{\tilde\alpha}w_{up}}^{^{\tilde\alpha}w_{down}}u|_{^{\tilde\alpha}S_{up}} dw.$$
The proof of {\it d)} therefore reduces to show that $-u_{rest}-(\tilde\alpha-\alpha)>0$ for all $\tilde\alpha\in[\alpha,\alpha+\Delta\alpha^*)$. It suffices to prove this for $\tilde\alpha=\alpha+\Delta\alpha^*$. By {\it b)}, $\Delta\alpha^*=-u_{rest}^3+\beta u_{rest}-\alpha$, so we have to prove that $u^3_{rest}-(\beta+1)u_{rest}+\alpha>0$, which was proved in Lemma \ref{LEM: sing hetero connection}{\it b)}.
\hfill$\square$\\

A straightforward corollary of Lemma \ref{LEM: sing hetero connection 2}, obtained by replacing $\tilde\alpha$ by $\alpha+z$, is that the singular ($\varepsilon_s=0$) phase portrait of (\ref{EQWcRd3DTravODE}) is as depicted in Figure \ref{FIG4}. Let $T$ denote the (2-dimensional) critical manifold of (\ref{EQWcRd3DTravODE}). Slices $^zT\subset T$ at $z\in[z_{rest}-\delta,z^*+\delta]$ are given by embedding in $\mathbb R^5$ the critical manifold of (\ref{sEQ: FHN wcusp ODE layer2}) for $\tilde\alpha=\alpha+z$. For each $z\in[z_{rest}-\delta,z^*+\delta]$ we denote with $^zT_{down/up}^{+/-}$ the embedding of $^{\tilde\alpha}S_{up/down}^{+/-}$, $\tilde\alpha=\alpha+z$, defined as in Figure \ref{sFIG: burst constr} and, similarly, with $^z\bar T_{down/up}^{+/-}$ the embedding of $^{\tilde\alpha}\tilde S_{up/down}^{+/-}$. For each $z\in[z_{rest}-\delta,z^*+\delta]$ we denote with $^z\mathcal H_{up/down}^{+/-}$ and $^z\mathcal H^0_{up}$ the embedding of $^{\tilde\alpha}\mathcal H_{up/down}^{+/-}$ and $^{\tilde\alpha}\mathcal H^0_{up}$, where they exist.

This singular phase portrait provides a skeleton for the application of the theorem in \cite[Section 4]{Jones1994}. In particular, we can build a singular homoclinic trajectory from the resting point to itself consisting of finitely many jumps along the family of heteroclinic orbits $^{z}\mathcal H_{up/down}^{+/-}$ of the layer dynamics of (\ref{EQWcRd3DTravODElayer}) connected slow trajectories of the reduced dynamics (\ref{EQWcRd3DTravODEreduced}). This trajectory was sketched in Figure \ref{FIG4}. Transversality properties of $^{\tilde\alpha}\mathcal H_{up/down}^{+/-}$ will translate into the needed transversality condition \cite[Eqs. (4.2) and (4.3)]{Jones1994}.

For the existence of the singular homoclinic trajectory, we invoke Lemma \ref{LEM: sing hetero connection 2}{\it d)}. { It implies that for $\tilde\varepsilon_{us}>0$ and as long as $z<z^*$, the $z$-coordinate of two successive jumps-down and of two successive jumps-up is strictly uniformly monotonically increasing by an amount bounded from below by $C\tilde\varepsilon_{us}$. After the last jump-down the trajectory might cross the hypersurface $\{w=-\lambda\}$ at $z=z^*$, which would break the transversality required by the exchange lemma. However, by smooth dependence of trajectories on the model parameters \cite[Theorem D.5]{Lee2013}, this behavior is not generic in $\tilde\varepsilon_{us}>0$, which explains the genericity condition in the statement of the theorem. The trajectory then converges toward the quasi-steady state of the slow-fast subsystem (\ref{EQWcRd3DTravODE}a-\ref{EQWcRd3DTravODE}c) and then to rest.}

We now verify the two transversality conditions [Jones and Koppel, Eqs. (4.2) and (4.3)]. Let
$$T^0:=\bigcup_{z\in[-\delta,\alpha_{PF}(\frac{1}{3})-\alpha]}(^z\bar T^{-}_{down}\ \cup\ \, ^z\bar T^{+}_{down})\cup \bigcup_{z\in(\alpha_{PF}(\frac{1}{3})-\alpha,\Delta\alpha^*+\delta]}\,^z\bar T_{down}$$
and
$$T^1:=\bigcup_{z\in[-\delta,\Delta\alpha^*+\delta]}^z\bar T_{up},$$
that is, $T^0$ and $T^1$ are compact, connected, normally hyperbolic (for the layer dynamics (\ref{EQWcRd3DTravODElayer})) submanifolds of the critical manifold $T$ that contains all the base and landing point of the heteroclinic orbits $^z\mathcal H_{up/down}^{+/-}$, in particular, $T^0$ contains base points of upward heteroclinic and landing points of downward heteroclinic, and viceversa for $T^1$.

To verify \cite[Eq. (4.2)]{Jones1994}, we have to prove that, in the layer dynamics (\ref{EQWcRd3DTravODElayer}), the first heteroclinic jump of the singular homoclinic trajectory, along $^{z_{rest}}\mathcal H_{up}^-$ is obtained as the transverse intersection $W^u_{rest}\cap_T W^s(T^1)$, where $W^u_{rest}$ is the 2-dimensional unstable manifold of the line of fixed points $\{(u_{rest},0,u_{rest},c,0),\ c\ {\rm near}\ c^*\}$, whereas $W^s(T^1)$ is 4-dimensional. Note that $W^u_{rest}\subset\{z=z_{rest}\}$ and that $W^s(T^1)\cap\{z=z_{rest}\}=W^s(S^{up})$, where $S^{up}$ is defined as in Lemma \ref{LEM: sing hetero connection}. It follows by Lemma \ref{LEM: sing hetero connection}{\it d)} that $W^u_{rest}\cap_T W^s(T^1)$ inside $\{z=z_{rest}\}$. Since the full space only adds the $z$-direction, $W^u_{rest}\cap_T W^s(T^1)$, which verifies  \cite[Eq. (4.2)]{Jones1994} for (\ref{EQWcRd3DTravODE}).

Next, we verify \cite[Eqs. (4.3)]{Jones1994} for all successive jumps. We write the condition explicitly for the first jump down along $^{\underline{z}}\mathcal H_{down}^+$ for some $\underline{z}>z_{rest}$. Computations for other jumps are similar and therefore omitted. We have to verify that $W^u(T^1|_{sing.orbit})\cap_T W^s(T^0)$. First, we note that the singular homoclininc trajectory is in the hypersurface $\{c=c^*\}$. Second, as $\tilde\varepsilon_{us}\to 0$, the slow motion after the first jump-up converges in the Hausdorff distance to the hypersurface $\{z=z_{rest}\}$, that is, $\lim_{\tilde\varepsilon_{us}\to0}T^1|_{sing.orbit}=T^1|_{c=c^*,z=z_{rest}}$, where the limit is again in the Hausdorff distance. Since $W^u(T^1|_{c=c^*,z=z_{rest}})\subset\{z=z_{rest}\}$, reasoning as above we can invoke  Lemma \ref{LEM: sing hetero connection}{\it d)} to conclude that $W^u(T^1|_{c=c^*,z=0})\cap_T W^s(T^0)$ inside the subspace $\{z=z_{rest}\}$ and therefore, by adding the $z$-dimension, in the whole space. Because $W^u(T^1|_{sing.orbit})$ is also $O(\tilde\varepsilon_{us})$-close to $W^u(T^1|_{c=c^*,z=z_{rest}})$, for $\tilde\varepsilon_{us}$ sufficiently small also $W^u(T^1|_{sing.orbit})\cap_T W^s(T^0)$.

We can now apply the theorem in \cite[Section 4]{Jones1994} to conclude that, for $\varepsilon_s>0$ sufficiently small, there exists a unique homoclinic solution of (\ref{EQWcRd3DTravODE}), which proves the statement of the theorem for $\beta=1/3$ and $\gamma=0$. For $\beta\neq 1/3$ and $\gamma\neq 0$, the statement follows from the persistence of transverse intersections to arbitrary $C^1$ perturbations.\hfill$\square$

\bibliographystyle{siamplain}
\bibliography{AF}

\end{document}